\documentclass[12pt]{article}
\usepackage{latexsym,amsfonts,amsmath,amssymb,amsthm}
\usepackage{amsthm,amsfonts,amsmath,amssymb}
\usepackage[utf8]{inputenc}
\usepackage[francais,english]{babel}
\usepackage{graphicx}
\usepackage{tikz}
\usepackage{nicefrac}
\usepackage{stmaryrd}
\usepackage{mathrsfs}
\usepackage{fourier}
\usepackage[all]{xy}
\usetikzlibrary{matrix,arrows}

\newtheorem{Proposition}{Proposition}[section]
\newtheorem{theorem}[Proposition]{Theorem}
\newtheorem{Definition}[Proposition]{Definition}
\newtheorem{Lemma}[Proposition]{Lemma}
\newtheorem{Remark}[Proposition]{Remark}
\newtheorem{Theorem}{Theorem}

\newtheorem{Corollary}[Proposition]{Corollary}

\renewcommand{\qed}{\nobreak \ifvmode \relax \else
      \ifdim\lastskip<1.5em \hskip-\lastskip
      \hskip1.5em plus0em minus0.5em \fi \nobreak
      \vrule height0.75em width0.5em depth0.25em\fi}
      
\newcommand\blfootnote[1]{%
  \begingroup
  \renewcommand\thefootnote{}\footnote{#1}%
  \addtocounter{footnote}{-1}%
  \endgroup
}

\newcommand{\TC}{ \tilde{\C}}

\newcommand{\TD}{ X_\Delta}

\newcommand{\ttor}{( \mathbb{C}^{\ast} )^{2}}

\newcommand{\C}{\mathcal{C}}
\newcommand{\A}{ \mathcal{A}}
\newcommand{\Ap}{ \A_\phi}
\newcommand{\T}{ \mathbb{T}}
\newcommand{\R}{ \mathbb{R}}

\newcommand{\rp}[1]{\mathbb{R} P^{#1}}
\newcommand{\F}{F}
\newcommand{\oF}{\overline{\F}}
\newcommand{\oC}{\overline{\C}}
\newcommand{\gp}{\gamma_\phi}

\newcommand{\Alp}{\Al_\phi}

\newcommand{\Cb}{\C_{\text{bl}}}
\newcommand{\Cbb}{\C_{B}}

\newcommand{\Co}{\C_{0}}
\newcommand{\sone}{(S^1)^2}
\newcommand{\phim}{\phi : \oC \rightarrow \TD}
\newcommand{\Argp}{ \Arg_\phi}
\newcommand{\Areap}{ \Area_\phi}
\newcommand{\CC}{\mathbb{C}}
\newcommand{\Z}{\mathbb{Z}}
\newcommand{\cp}[1]{\CC P^{#1}}
\newcommand{\D}{\mathcal{D}}
\newcommand{\XD}{X_\Delta}
\newcommand{\cZ}{\mathcal{Z}}
\renewcommand{\O}{\mathscr{O}}
\newcommand{\VS}{VN\:}

\DeclareMathOperator{\itr}{int}
\DeclareMathOperator{\dif}{d}
\DeclareMathOperator{\Log}{Log}

\DeclareMathOperator{\Arg}{Arg}
\DeclareMathOperator{\Al}{Alga}
\DeclareMathOperator{\Area}{Area}
\DeclareMathOperator{\grad}{grad}
\DeclareMathOperator{\ord}{ord}
\DeclareMathOperator{\New}{New}
\DeclareMathOperator{\conv}{conv}
\DeclareMathOperator{\Sub}{Sub}
\DeclareMathOperator{\Abs}{Abs}

\DeclareMathOperator{\ind}{ind}
\DeclareMathOperator{\Top}{Top}
\newcommand{\vv}{\nu}

\newenvironment{psmallmatrix}
  {\left(\begin{smallmatrix}}
  {\end{smallmatrix}\right)}

\newcommand{\nocontentsline}[3]{}
\newcommand{\tocless}[2]{\bgroup\let\addcontentsline=\nocontentsline#1{#2}\egroup}

\title{A generalisation of simple Harnack curves}
\author{Lionel Lang}

\begin{document}

\maketitle

\begin{abstract}
In this paper, we suggest the following generalisation of Mikhalkin's simple Harnack curves: a generalised simple Harnack curve is a parametrised real algebraic curve in $\ttor$ with  totally real logarithmic Gauss map. We investigate which of the many properties of simple Harnack curves survive the latter generalisation. We also show how tropical geometry allows to construct plenty of examples. Since generalised Harnack curves can develop arbitrary singularities, in contrast with the original definition where only real isolated double points can appear,  we pay a special attention to the simplest new instance of generalised Harnack curves, namely curves with a single hyperbolic node. In particular, we give their topological classification as in \cite{Mikh} and show how such curves can be recovered from their spine.
\end{abstract}

\blfootnote{During the writing of this paper, the author was supported by the FNS project 140666, the ERC grant TROPGEO and SwissMAP. The author is grateful to  G.Mikhalkin, R.Crétois  and E.Brugallé for many useful discussions and comments.}

\begin{large}
\textbf{Introduction}
\end{large}
\\

Simple Harnack curves introduced in \cite{Mikh} are at the crossroad between real algebraic geometry, symplectic geometry, complex analysis and physics, see for instance \cite{Mikh}, \cite{Bru14}, \cite{Kri} and \cite{KOS}. Among the many characterisations of simple Harnack curves, the most concise one is given in term of the amoeba map $\A:\ttor \rightarrow \R^2$ defined by 
\[\A(z,w):=\big(\log\vert z\vert, \log\vert w\vert\big).\]
A real algebraic curve $A\subset \ttor$ is a simple Harnack curve if $\A_{\vert A}$ is at most $2$-to-$1$, see \cite{MR}. The original definition of \cite{Mikh} is topological and involves the natural toric compactification $\ttor \subset \XD$ given by ``the'' Newton polygon $\Delta$ of $A$. Strikingly enough, the real curve $A$ is a simple Harnack curve if and only if it is maximal in any of the following ways:

-- The curve $A$ is maximal with respect to a finite collection of Smith-Thom inequalities, see \cite{Mikh04}.

-- The area of $\A(A)$ is maximal among curves with the same Newton polygon, see \cite{MR}.

-- The total logarithmic curvature of $\R A$ is maximal among curves with the same Newton polygon, see \cite{PR}.\\
An other useful characterisation of simple Harnack curves is given in \cite{MO}, see also Section \ref{secdef} below. Originally, the simple Harnack curves introduced in \cite{Mikh} were assumed to be smooth. Later on, the authors of \cite{MR} introduced singular simple Harnack curves as deformation of the curves of \cite{Mikh}.
In \cite{KO}, a complete description of the space of simple Harnack curves in $\cp{2}$ were given (see \cite{Ola} for a generalisation to any toric surface). It turned out that the singular simple Harnack curves of \cite{MR} are exactly the curves at the boundary of the space of smooth simple Harnack curves. In particular, simple Harnack curves can only develop real isolated double points as singularities.

In the present paper, we suggest the following generalisation: A real algebraic map $\phi: \C \rightarrow \ttor$ from a real smooth punctured algebraic curve $\C$ is a generalised Harnack curve if the induced logarithmic Gauss map $\gp: \C \rightarrow \cp{1}$ is totally real, see  Section \ref{secdef}. Recall that for an algebraic curve $A \subset \ttor$ given as the vanishing locus of a polynomial $f$, the logarithmic Gauss map $\gamma: A \dashrightarrow \cp{1}$ defined on the smooth locus of $A$ is given by 
\[\gamma(z,w)= \big[ z \cdot \partial_z f(z,w)\, ;\;  w \cdot \partial_w f(z,w)\big]. \]
The parametrisation $\phi$ allows to extend the map $\gamma$ to a meromorphic function $\gp$ on $\C$. The fact that the above definition generalises the notion of simple Harnack curve of \cite{Mikh} and \cite{MR} follows from the works \cite{Mikh} and \cite{PR}.

In this text, we investigate the properties of generalised Harnack curves. In Section \ref{secdef}, we show that for any generalised Harnack curve $\phi: \C \rightarrow \ttor$, the real part $\R \oC$ of the compactified curve $\oC$ is an M-curve, see Theorem \ref{thm:max}. We also show that the characterization \cite[Theorem A1]{MO} extends to generalised Harnack curves, see Theorem \ref{thmalga}. In Section \ref{sec:trop}, we provide a tool for constructing generalised Harnack curves. To do so, we introduce the notion of tropical Harnack curves. We show in Theorem \ref{thmexistharncurv} that the approximation of the latter tropical curves leads to generalised Harnack curves. Roughly speaking, this construction is a parametrised version of the combinatorial patchworking with no twisted edges, see for instance \cite[\S 3]{BIMS}. The latter construction allows in particular to produce generalised Harnack curves with hyperbolic nodes and complex conjugated nodes. In Section \ref{sec:1hn}, we focus on generalised Harnack curves with a single hyperbolic node. As noticed earlier, the latter curves can not be obtained as deformation of simple Harnack curves. As in \cite{Mikh}, we undertake the classification of the topological pairs
\[\Top(\phi):=\Big( \R \TD, \; \R \C  \bigcup_{j\in \Z/n\Z} \R \mathcal{D}_j \Big)  \]
where $\phi:\C:\rightarrow \ttor$ is a generalised Harnack curves with a single hyperbolic node, $\Delta$ is ``the'' Newton polygon of $\phi(\C)$ and the $\D_j$ are the toric divisors of the toric surface $\XD$.
In Theorem \ref{thmtopclass}, we show that the possible topological pairs $\Top(\phi)$ are indexed by the vertices of $\Delta$ corresponding to smooth points of $\XD$. The latter theorem is obtained as a consequence of Theorem \ref{thmgoodspine} stating that the spine of such a generalised Harnack curve is always a tropical Harnack curve with a single hyperbolic node, see Definitions \ref{def:tropedgeleafnode} and \ref{deftropsimpleHarnack}.

To conclude this introduction, let us advertise the class of generalised Harnack curves introduced in this paper. To underline how important the latter class of curve might be, observe that generalised Harnack curves encompass Mikhalkin's simple Harnack curves as well as reduced A-disriminant curves. Indeed, the logarithmic Gauss map of A-disriminant curves is a birational isomorphism according to the Horn-Kapranov parametrisation, see for instance \cite{PasNil}. For this reason, generalised Harnack curves might be thought of as an interpolation between the two aforementioned classes of curves. In this paper, we establish the foundation to investigate further connections between generalised Harnack curves and their tropical avatars. For instance, it would be very useful to generalise the local coordinates \cite[Propositions 6 and 10]{KO} in order to study the deformations of generalised Harnack curves. In particular, it would have interesting applications to the study of Severi varieties, in the fashion of \cite{CL1} and \cite{L19_2}. Finally, let us mention that the definition of simple Harnack curves has been generalised to higher dimension, see for instance \cite{Mikh04}. Such generalisation were studied in  \cite{BMRS} where the authors showed that smooth hypersufaces in $(\mathbb{C}^*)^n$ with totally real logarithmic Gauss map are limited to hyperplanes. As suggested by the present work, one should investigate the case of singular hypersurfaces, considering the fact that reduced A-discriminant hypersurfaces already provide a non-trivial example of hypersurfaces with (almost) totally real logarithmic Gauss map.

\tableofcontents

\section{Setting}

\subsection{Curves in toric surfaces}\label{sec:loggeo}

Consider the complex torus $\ttor$ with coordinates $(z,w)$. In this text, the symbol $\Delta$ refers to a \textbf{compact lattice polygon} $\Delta\subset \R^2$ with non empty interior given as the convex hull of a finite set of points in the character lattice $\Z^2$ of the torus $\ttor$. Define the \textbf{toric surface} $X_\Delta \supset \ttor$ as the closure of the embedding 
$$\ttor \hookrightarrow \cp{\vert \Delta \cap \Z^2 \vert-1}$$
given coordinate-wise by the monomials $z^\alpha w^\beta$ for $(\alpha, \beta) \in \Delta \cap \Z^2$, see \cite[\S 2.3]{CLS}. We denote by $n$ the \textbf{number of sides of} $\Delta\subset \R^2$ and fix once and for all a cyclic indexation in $\Z/n\Z$ of these sides by going along $\partial \Delta$ counter-clockwise. For any $j\in \Z/ n\Z$, we denote by $\Delta_j$ the $j^{th}$ \textbf{side of } $\Delta$ and by $\D_j$ the \textbf{toric divisor} corresponding to $\Delta_j$. Recall that $\XD$ is smooth at the point $\D_j\cap\D_{j+1}$ if and only if the primitive vectors supporting $\Delta_j$ and $\Delta_{j+1}$ generate $\Z^2$, see \cite[\S 2.5]{Ful}. In such case, we will say that the vertex $\Delta_j\cap\Delta_{j+1}$ of $\Delta$ is \textbf{smooth}. We will also denote 
\[ b_\Delta := \# (\partial \Delta \cap \mathbb{Z}^2) \quad \text{and} \quad g_\Delta := \# (\itr(\Delta) \cap \mathbb{Z}^2). \]
The complex conjugation induced by the coordinates $(z,w)$ on $\ttor$ extends to a complex conjugation on $X_\Delta$. For a complete, smooth real algebraic curve $\oC$, we use $\phi : \oC \rightarrow \TD$ to refer to a \textbf{unramified real algebraic map} and denote $\C := \phi^{-1}\big(\ttor\big)\subset\oC$. We will say that $\phi$ has \textbf{degree} $\Delta$ if $\phi(\C)$ is the vanishing locus of a Laurent polynomial with Newton polygon $\Delta$. In particular, the curve $\phi(\oC) \subset \XD$ intersects every divisor $\mathcal{D}_j$ in   
$\#(\Delta_j \cap \mathbb{Z}^2) -1$ many points, counted with multiplicities. Also, the curve $\phi(\oC)$ does not contain any toric fixed point of $\XD$. For practical reasons, let us introduce the following notations 
$$ \C_{j} := \phi^{-1}(\D_j) \subset \oC \text{ for any } j \in \Z/n\Z \; \text{ and } \; \C_\infty := \bigcup_{j \in \Z/n\Z} \C_j.$$
The integer $b_\Delta$ defined above is then the intersection multiplicity of $\phi(\oC)$ with the union of all the toric divisors. In particular, if $\phi(\oC)$ intersects each of them transversally, then $ \#\C_\infty= b_\Delta$.
According to \cite{Kho}, the integer $g_\Delta$ is the arithmetic genus of $\phi(\oC) \subset \XD$. 

Recall the \textbf{amoeba map} 
\[
\begin{array}{rcl}
\A \: : \: \ttor & \rightarrow & \mathbb{R}^2 \\
(z,w) & \mapsto & \big( \log \vert z \vert, \log \vert w \vert \big)
\end{array}
\]
and the coordinate-wise argument map
\[
\begin{array}{rcl}
\Arg \: : \: \ttor & \rightarrow & (S^1)^2 \\
(z,w) & \mapsto & \big( \arg(z), \arg(w) \big)
\end{array}.
\]
Denote respectively $\Ap := \A \circ \phi$ and $\Argp := \Arg \circ \, \phi$ the maps defined on $\C$.
\begin{Definition}
The \textbf{amoeba} of \,$\C$ (relative to $\phi$) is the subset $ \Ap (\C) \subset \mathbb{R}^2$. The \textbf{coamoeba} of $\C$ (relative to $\phi$) is the subset $ \Argp ( \C) \subset \big(S^1\big)^2$.
\end{Definition}

\begin{Proposition}[see \cite{FPT}]
The amoeba of $\C$ is a closed subset of $\;\mathbb{R}^2$ such that  every connected component of $\; \mathbb{R}^2\setminus\Ap(\C)$ is convex.
\end{Proposition}

The coordinates $z$ and $w$, seen as meromorphic functions on $\oC$, induce respectively the two meromorphic forms $d \log z$ and $d \log w$ on $\oC$. Define the \textbf{logarithmic Gauss map} 
\[
\begin{array}{rcl}
\gp : \oC & \rightarrow & \cp{1} \\
p &\mapsto& \left[ - \big(\dif\log w\big)(p) \, ; \; \big(\dif\log z\big) (p) \right]
\end{array}
\]
as the quotient of these two forms. The latter definition agrees with the original definition of \cite{Kap} when $\phi$ is an embedding.

\begin{Proposition}\label{degloggauss}
The degree of the logarithmic Gauss map $\gp$ is at most $ - \chi (\C)$ with equality if and only if the map $\phi$ is an immersion. 
\end{Proposition}
\noindent
\begin{proof}
The degree of the logarithmic Gauss map $\gp$ is equal to the number of points in $\gp^{-1}\big( [u;v] \big)$ for a generic point $[u;v] \in \cp{1}$. In turn, the set $\gp^{-1}\big( [u;v] \big)$ is a subset of the set of zeroes of the form $ u \cdot \dif\log w + v \cdot \dif\log z $. For a generic $[u;v] \in \cp{1}$, the zeroes of $ u \cdot \dif\log w + v \cdot \dif\log z $ that are not in $\gp^{-1}\big( [u;v] \big)$ are exactly the common zeroes of  $\dif\log w=\frac{w'}{w} dw$ and  $\dif\log z=\frac{z'}{z} dz $, that is the set of points where $\phi$ is not an immersion. In order to prove the statement, it remains to prove that the number of zeroes of  $ u \cdot \dif\log w + v \cdot \dif\log z $ is exactly $- \chi (\C)$. To see this, note on the one hand that the form $ u \cdot \dif\log w + v \cdot \dif\log z $ has no pole on $\C$. On the other hand, every point in $\C_\infty$ is a simple pole of this form whenever $[u;v]$ is not in $\rp{1}$. Indeed, both $ \dif\log w$ and  $\dif\log z $ have a simple pole at each such point with  residue in  $2\pi i \Z$. Now, we know by Riemann-Roch that the divisor of a meromorphic form on $\C$ has degree $g-2$ where $g$ is the genus of $\C$. It follows that the form $ u \cdot \dif\log w + v \cdot \dif\log z $ has $g-2+ \# \C_\infty$ many zeroes. Since $- \chi (\C)= g-2+ \# \C_\infty$, the result follows. 
\end{proof}

The maps $\Ap$ and $\Argp$  are differentiable maps between smooth surfaces. Following \cite{Mikh}, we denote by $\F \subset \C$ the closure of the set of points where the tangent map $T \Ap$ has rank $1$. Observe that the points where $\phi$ is not an immersion need not to be in $F$.

\begin{Lemma}\label{lem:critloc}
The tangent map $T \Ap$ has rank $1$ at $p$ if and only if the map $T \Argp$ has rank $1$ at $p$. In particular, the subset $\F \subset \C$ is the closure of the set of points where $T \Argp$ has rank $1$. Moreover, we have
\[ \overline{F} = \gp^{-1}\big(\rp{1} \big).\]
\end{Lemma}

\begin{proof}
We borrow the arguments of \cite[Lemma 3]{Mikh}. Let $p \in \C$  be a point such that $T \Ap$ has rank $1$ at $p$. In particular, the map $\phi$ is an immersion at $p$. Then, the tangent space  $T_{\phi(p)} \phi(\C)$ contains a vector $v\in\mathbb{C}^2$ tangent to the torus $\big\lbrace (z,w) \in \ttor \; \big| \; \vert z \vert = \vert p_1 \vert, \; \vert w \vert = \vert p_2 \vert\big\rbrace$. Equivalently, the vector $v$ has purely imaginary logarithmic coordinates. Therefore, the vector $i\cdot v \in T_{\phi(p)} \phi(\C)$ has real logarithmic coordinates. Thus, the vector $i\cdot v$ is tangent to $\big\lbrace (z,w) \in \ttor \; \big| \; \arg(z) = \arg(p_1), \; \arg(w)= \arg(p_2) \big\rbrace$. Equivalently, the vector $i\cdot v$ is in the kernel of $T \Argp$. We deduce that $T \Argp$ has rank $1$ at $p$ since $v$ is not in the kernel of $T \Argp$. The latter arguments are obviously symmetric so that $T \Ap$ has rank $1$ at $p$ as soon as $T \Argp$ as rank $1$ at $p$. The first part of the statement follows. 
We know that $T \Ap$ has rank $1$ at $p$ if and only if $\phi$ is an immersion at $p$ and  $T_{\phi(p)} \phi(\C)$ contains a vector $i\cdot v$ with real coordinates, if and only if $\phi$ is an immersion at $p$ and $\gp(p) \in \rp{1}$. Since the set of points where $\phi$ is not an immersion is discrete, the second part of the statement follows.
\end{proof}

\begin{Remark}
Any point of $\C$ mapped to the boundary of $ \Ap(\C)$ belongs to $F$. By the above lemma, 
the real part $\mathbb{R} \oC$ of the curve $\oC$ is always contained in  $\overline{F}$.
\end{Remark}

Write $z=e^{x_1+i x_2}$ and $w=e^{y_1+i y_2}$. The following observation is due to Mikhalkin.

\begin{Lemma}
For any algebraic map $\phim$, the 2-forms $\Ap {^\ast}\big(dx_1 \wedge dy_1\big)$ and $\Argp^{\;\,\ast}\big( dx_2 \wedge dy_2 \big)$ coincide on $\C$. It implies that
\[  \int_\C \Ap {^\ast}\big(dx_1 \wedge dy_1\big) = \int_\C \Argp^{\;\,\ast}\big(dx_2 \wedge dy_2\big) = 0  \]
and
\[  \int_\C \left| \Ap {^\ast}\big(dx_1 \wedge dy_1\big) \right| = \int_\C \left| \Argp^{\;\,\ast}\big(dx_2 \wedge dy_2\big) \right|. \]
\end{Lemma}
\noindent
\begin{proof}
At any point of $\C$, we have the local descriptions $\Ap = \Re \circ \Log \circ \phi$ and $\Argp = \Im \circ \Log \circ \phi$ where $\Log$ is a local determination of the coordinatewise logarithm and $\Re$ and $\Im$ are the projections on the real and imaginary coordinates respectively. Set $x:=x_1+i x_2$ and $y:=y_1+i y_2$. Then, the pullback to the curve $\C$ of the holomorphic $2$-form $dx\wedge dy$ by $\Log \circ \phi$ vanishes everywhere. In particular, we have that the pullback of $\Re\big(dx\wedge dy\big)=dx_1 \wedge dy_1-dx_2 \wedge dy_2$ vanishes everywhere on $\C$. The first part or the statement is proven. The equalities $\int_\C \Ap {^\ast}\big(dx_1 \wedge dy_1\big) = \int_\C \Argp^{\;\,\ast}\big(dx_2 \wedge dy_2\big)$ and $\int_\C \left| \Ap {^\ast}\big(dx_1 \wedge dy_1\big) \right| = \int_\C \left| \Argp^{\;\,\ast}\big(dx_2 \wedge dy_2\big) \right|$ follow directly. Finally, the smooth map $\Ap : \C \rightarrow \R^2$ has a well defined degree as the map is proper. This degree is necessarily $0$ as $\Ap$ is not surjective. It implies that $\int_\C \Ap {^\ast}\big(dx_1 \wedge dy_1\big)= \deg (\Ap) \Area\big( \Ap(\C)\big)=0$ and the lemma is proven.
\end{proof}

\begin{Definition}\label{def:areap}
For an algebraic map $\phim$, define 
\[\Areap (\C) := \int_\C \left| \Ap {^\ast}\big(dx_1 \wedge dy_1\big) \right| = \int_\C \left| \Argp^{\;\,\ast}\big(dx_2 \wedge dy_2\big) \right|. \]
\end{Definition}

Define the \textbf{Ronkin function} $N_\phi : \R^2 \rightarrow \R$ associated to $\phi$ by 
\[ N_\phi (x,y)  :=  \frac{1}{(2i\pi)^2} \displaystyle \int_{\mathcal{A}^{-1}(x,y)} \dfrac{\log \vert f(z,w) \vert }{zw} dz \wedge dw \]
where $f(z,w)$ is any Laurent polynomial with Newton polygon $\Delta$  and such that $\phi(\C)= \big\lbrace (z,w) \in \ttor \, \big| \linebreak \;  f(z,w)=0 \big\rbrace$. The function $N_\phi$ is convex and affine linear on each connected component of $ \R^2 \setminus \Ap(\C)$, see \cite[Theorem p. 483]{PR}. Consequently, the gradient $\grad \, N_\phi $ is constant on each such component. Define the \textbf{order map} $\ord : \big\lbrace \text{connected components of } \R^2 \setminus \Ap(\C)\big\rbrace \rightarrow \Z^2$ by 
\[\ord (E):= \grad N_\phi(p)\]
for any point $p\in E$. The latter definition agrees with the original definition \cite[Definition 2.1]{FPT} according to \cite[Theorem p. 483]{PR}. The proposition below is a follows from \cite[Propositions 2.4 and 2.5]{FPT}.

\begin{Proposition}
The order map $\ord$ is injective and valued in $\Delta \cap \Z^2$. A components $E \subset \mathbb{R}^2 \setminus \mathcal{A}(\C) $ maps to $ \itr(\Delta) \cap \mathbb{Z}^2$ if and only if $E$ is compact. Moreover, every vertex of $\Delta$ is in the image of $\ord$.
\end{Proposition}

Denote by $\mathscr{E}_\phi$ the \textbf{set of connected components of} $ \R^2 \setminus \Ap(\C)$. For any  $E \in \mathscr{E}_\phi$, denote by $N_\phi^E$ the affine linear function on $\mathbb{R}^2$ extending $(N_\phi)_{\vert_E}$. Then, the \textbf{spine} $\mathscr{S}_\phi$ is defined as the corner locus of the piecewise affine linear and convex function
\[S_\phi := \max_{E \in \mathscr{E}_\phi} N_\phi^E.\]
The spine $\mathscr{S}_\phi \subset \R^2$ is a piecewise linear graph in the plane. In Section \ref{sec:trop}, we will see that $\mathscr{S}_\phi$ is a tropical curve when equipped with appropriate weights. According to \cite[Theorem 1]{PR}, we have the following.

\begin{theorem}\label{thmPR}
The spine $\mathscr{S}_\phi$ is a deformation retract of the amoeba $\Ap(\C)$.
\end{theorem}

Let us conclude this section with some elementary lemmas.

\begin{Lemma} \label{valinf}
Let $(a,b) \in \mathbb{Z}^2$ be a primitive integer vector supporting the side $\Delta_j$. For any point $p \in \C_j$, we have
\[ \gp (p) = \left[a \, ; \; b\right]. \]
In particular, we have the inclusion $\C_\infty \subset \oF$.
\end{Lemma}

\begin{proof}
Assume for a moment that $\Delta_{j-1}$ and $\Delta_j$ are supported respectively on the on the vertical and on the horizontal coordinate axis of $\Z^2$. In particular, we have that $\Delta_{j-1}\cap\Delta_j=(0,0)$ and that $\ttor \subset \CC^2$ is an affine chart of $\XD$ in a neighbourhood of $\D_{j-1}\cap \D_j$ such that $\D_j=\{z=0\}$ and $\D_{j-1}=\{w=0\}$. Then, there exists a holomorphic coordinate $t$ on $\oC$ centered at $p$ such that 
\[ z(t)=z_0+z_n t^n + o(t^n) \; \; \text{ and } \; \;  w(t)=t^m + o(t^m)\] 
where $z_0, z_n \in \CC^*$ and $n, m \in \Z_{\geq 1}$. It follows that 
\[
\begin{array}{rl}
\gp (p) & = \; \displaystyle \lim_{t \rightarrow 0} \left[ - \dif \log \big(w(t)\big) \,
; \;  \dif \log \big(z(t)\big) \right] \; = \; \displaystyle \lim_{t \rightarrow 0} \left[ - \dfrac{w'(t)}{ w(t)} \, ; \;  \dfrac{ z'(t) }{z(t)} \right]  \\
					&  \\
					& = \;  \displaystyle \lim_{t \rightarrow 0} \left[ - \dfrac{mt^{m-1}+o(t^{m-1})}{t^m + o(t^m)} \, ; \;  \dfrac{nz_n t^{n-1} + o(t^{n-1})}{z_0+z_n t^n + o(t^n)} \right] \; = \; \displaystyle \lim_{t \rightarrow 0} \left[ - \dfrac{mt^m+o(t^{m})}{t^m + o(t^m)} \, ; \;  \dfrac{nz_n t^{n} + o(t^{n})}{z_0+z_n t^n + o(t^n)} \right] \\
					\\
					&= \; \displaystyle \lim_{t \rightarrow 0} \left[ - \dfrac{m+o(1)}{1 + o(1)} \, ; \;  \dfrac{nz_n t^{n} + o(t^{n})}{z_0+z_n t^n + o(t^n)} \right] \; = \; \left[1\, ; \; 0\right]. 
\end{array}
\]

In the general case, consider the integer affine linear transformation $(\alpha,\beta)\mapsto (\alpha,\beta) \cdot \begin{psmallmatrix} d & -b \\ -c & a \end{psmallmatrix} + (e,f)$ of the character lattice $\Z^2$ sending the polygon $\Delta$ to a polygon $\Delta'$ such that the respective images of $\Delta_{j-1}$ and $\Delta_j$ are as in the above paragraph. Here, we have $\det \begin{psmallmatrix} d & -b \\ -c & a \end{psmallmatrix}=1$. The dual map of complex tori $\psi : \ttor \rightarrow \ttor$ given by $\psi(z,w)=\big(z^dw^{-c}, z^{-b}w^a\big)$ extends to an isomorphism $\psi: X_{\Delta'} \rightarrow \XD$. In particular, we have a factorisation $\phi= \psi \circ \phi'$ where $\phi' : \oC \rightarrow X_{\Delta'}$. It follows that 
\[
\begin{array}{rl}
\gp & = \; \displaystyle \left[ - \dif \log \big(z^{-b}w^a\big) \,
; \;  \dif \log \big(z^dw^{-c}\big) \right] \; = \; \displaystyle \left[ - \dif \log w \, ; \;  \dif \log z \right] \cdot \begin{pmatrix} a & c \\ b & d \end{pmatrix} \\
					&  \\
					& = \; \gamma_{\phi'} \cdot  \begin{pmatrix} a & c \\ b & d \end{pmatrix}.
\end{array}
\]
The result follows now from the above paragraph.
\end{proof}

\begin{Lemma}\label{lem1}
--  $ \oF \subset \oC $ is smooth if and only if $\gp$ has no branching point on $ \mathbb{RP}^1$. In this case, $\oF$ is a disjoint union of smoothly embedded circle in $\oC$.

-- If $\alpha$ is a smooth connected component of $\F$ and $(\Ap)_{\vert \alpha}$ is non-constant, then the Gauss map of the parametrised curve $\Ap : \alpha\rightarrow \R^2$ is given by the restriction of $\gp$ to $\alpha$. In particular, the latter Gauss map is monotonic and $(\Ap)_{\vert \alpha}$ has no inflection point.

-- If $\alpha$ is a smooth connected component of $\F$ and $(\Argp)_{\vert \alpha}$ is non-constant, then Gauss map of the parametrised curve $\Argp : \alpha \rightarrow \R^2$ is given by the restriction of $\gp$ to $\alpha$. In particular, the latter Gauss map is monotonic and $(\Argp)_{\vert \alpha}$ has no inflection point.
\end{Lemma}

\begin{proof}
The first part of the Lemma is proven in \cite[Proposition 1.1]{L19}. For the second and third parts, observe first that the restriction of $\gp$ to $\alpha$ is a local diffeomorphism, otherwise $\alpha$ would not be smooth. If the Gauss maps of  $(\Ap)_{\vert\alpha}$ and $(\Argp)_{\vert \alpha}$ are described by $\gp$ as we claim, it follows that the latter parametrised curves have no inflection. Note also that $(\Ap)_{\vert\alpha}$ and $(\Argp)_{\vert\alpha}$ are real analytic so that each map is either constant or locally injective. In particular, each parametrised curve admits a Gauss map under the present assumptions. According to the proof of Lemma \ref{lem:critloc},  the tangent line to $\oC$ at a point $p \in \alpha$ contains to real lines $\ker T_p \Ap$ and $\ker T_p \Argp$ obtained one from another by multiplication by $i$.  The image $v$ of $\ker T_p \Argp$ by the tangent map to $\Log \circ \phi$ (where $\Log$ is any branch of the coordinate-wise complex logarithm) is a vector with real coordinates in $T_{Log(\phi(p))}\CC^2$. In the latter logarithmic coordinates, the logarithmic Gauss map $\gp$ is given by the usual Gauss map (see \cite[\S 3.2]{Mikh}) while the map $\Ap$ is given by the projection onto the first factor of $\CC^2= \R^2 \oplus i \R^2$. Therefore, we have that $\gp(p)= [v] \in \cp{1}$. As $v$ has real coordinates, the tangent map to $\Ap$ is the identity on $v$ in the present logarithmic coordinates. It follows that the Gauss map to the parametrised curve $(\Ap)_{\vert \alpha}$ equals also $[v] \in \cp{1}$ at $p$. The second part of the lemma follows. The third part is proven similarly.
\end{proof}

Define $\C_{\Arg}$ to be the \textbf{real oriented blow-up} of $\oC$ at every point of $\C_\infty$. Denote by $S_p \subset \C_{\Arg}$ the fiber of the blow-up over $ p\in \C_\infty$. In particular, we have $S_p\simeq S^1$.

\begin{Lemma}\label{extarg}
The map $\Argp$ extends to $\C_{\Arg}$. Moreover, if the edge $\Delta_j$ is supported by a primitive integer vector $(a,b)$ and $\phi(p)$ belongs to $\mathcal{D}_j$, then the restriction $(\Argp)_{S_p}$ is an $m$-covering of a geodesic of slope $(-b,a)$, where $m$ is the intersection multiplicity of $\phi(\oC) \cap \D_j$ at $\phi(p)$.
\end{Lemma}

\begin{proof}
Assume first that $\Delta_{j-1}$ and $\Delta_j$ are supported respectively  on the vertical and on the horizontal coordinate axis of $\Z^2$. From the proof of Lemma \ref{valinf}, we know that there exists a local coordinate $t$ centered at $p$ such that 
\[z(t)=z_0+z_n t^n + o(t^n) \; \; \text{ and } \; \;  w(t)=t^m + o(t^m)\]
where $z_0, z_n \in \CC^*$ and $n, m \in \Z_{\geq 1}$. Here, the integer $m$ is the intersection multiplicity of $\oC \cap \D_j$ at $\phi(p)$.
For $t=r e^{i \theta}$, we have
\[z(t)=z_0+z_n r^n e^{i n\theta} + o(r^n) \; \; \text{ and } \; \;  w(t)=r^m e^{i m\theta} + o(r^m)\]
It follows that for any $e^{i\theta} \in S^1$
\[ \lim_{r \rightarrow 0} \arg \big(z( re^{i\theta}) \big)= \arg(z_0)
\; \; \text{ and } \; \; \lim_{r \rightarrow 0} \arg \big( w( re^{i\theta}) \big) =  e^{i m\theta}. \]
This proves the lemma in this case. For the general case, it suffices to use a change of coordinates $\psi$ as in the proof of Lemma \ref{valinf}.
\end{proof}

\subsection{Simple  Harnack curves}\label{sec:defshc}

In this section, we give  a brief account on simple Harnack curves as introduced in \cite{Mikh}. As we will mostly deal with singular curves in this paper, we adopt the parametrised point of view used in \cite{Bru14}.
Recall that the smooth real algebraic curve $\oC$ of genus $g$ is an M-curve if its real part $\mathbb{R} \oC \subset \oC$ has the maximal number of connected components, that is $g+1$.

\begin{Definition}{\cite[Definition 2]{Mikh}}\label{defmikh1}
The real algebraic map $\phi: \oC \rightarrow X_\Delta$ is a \textbf{simple Harnack curve} if
\begin{itemize}
\item[$a)$] $\phi$ is an embedding of degree $\Delta$ and $\oC$ is an M-curve,
\item[$b)$] there exists a connected component $\O \subset \R \oC$ such that $\C_\infty \subset \O$,
\item[$c)$] we can orient the topological circle $\O$ such that for any $j \in \Z/n\Z$, any $p_{j-1} \in \C_{j-1}$, $p_{j} \in \C_{j}$ and $p_{j+1} \in \C_{j+1}$, the point $p_{j}$ lies inside the open arc in $\O$ going from $p_{j-1}$ to $p_{j+1}$.
\end{itemize}
\end{Definition}

Observe that we do not require that $\phi(\oC)$ intersects $\D_j$ transversally for any $j \in \Z / n\Z$ in contrast with \cite[Definition 2]{Mikh}.

\begin{figure}[h]
\begin{center}
\input{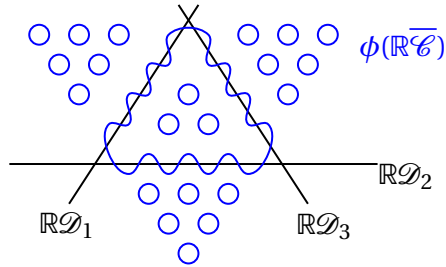}
\end{center}
\caption{A simple  Harnack curve of degree $8$ in $\rp{2}$.}
\label{Hk1}
\end{figure}

\begin{theorem}{\cite[Theorem 3]{Mikh}}
For any simple Harnack curve $\phi$ in $\TD$ intersecting every $\D_j$ transversally, the topological type  $$\Big( \R \TD\, , \; \, \R \phi\big(\oC\big) \bigcup_{j\in \Z/n\Z} \R \mathcal{D}_j \Big)$$ depends only on $\Delta$, see for instance Figure \ref{Hk1}.
\end{theorem}

\begin{Definition}{\cite[Definition 3]{MR}}\label{defmikh2}
The map $\phi: \oC \rightarrow X_\Delta$ is a \textbf{singular simple Harnack curve} if
\begin{itemize}
\item[$a)$] the only singular points of $\phi\big(\oC\big)$ are real isolated double points,
\item[$b)$] the operation of replacing of the singular points of $\phi\big(\oC\big)$ with small real ovals gives a simple Harnack curve in $\TD$.
\end{itemize}
\end{Definition}

\begin{theorem}{\cite[Theorem 1]{MR}}
The real algebraic map $\phi: \oC \rightarrow X_\Delta$ of degree $\Delta$  is a possibly singular simple Harnack curves if and only if one of the following conditions is satisfied 
\begin{itemize}
\item[$a)$] $\Ap : \C \rightarrow \mathbb{R}^2$ is at most 2-to-1.
\item[$b)$] $ \Area \big( \Ap (\C) \big)\big) = \pi^2 \Area (\Delta)$.
\end{itemize}
\end{theorem}

\begin{figure}[h]
\begin{center}
\input{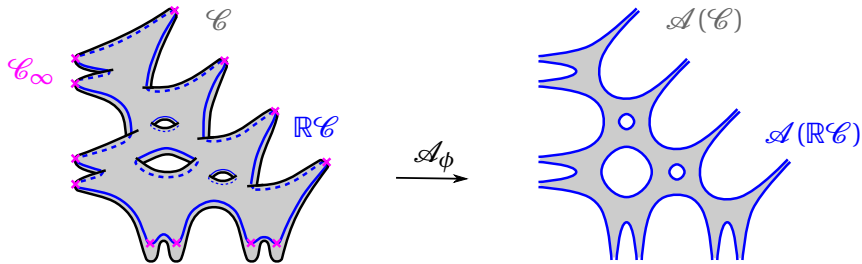}
\end{center}
\caption{The amoeba map on a simple  Harnack quartic.}
\label{Amap1}
\end{figure}

Observe that the above theorem is accurate only when (possibly singular) simple Harnack curves are allowed to intersect the divisors $\D_j$ non-transversally. The property $b)$ is illustrated in Figure \ref{Amap1}. There is yet an other equivalent characterisation of simple Harnack curves that is of special interest to us.

\begin{theorem}[see \cite{Mikh} and \cite{PRis}]\label{thmPR}
A real algebraic embedding $\phi : \oC \rightarrow \TD$ is a simple Harnack curves if and only if the logarithmic Gauss map $ \gp : \oC \rightarrow \cp{1}$ is totally real, that is 
\[ \gp^{-1}\big(\rp{1}\big) = \R \oC. \]
\end{theorem}

\section{Generalisation and first properties}\label{secdef}

\begin{Definition}\label{defsimpleHarnack}
A real algebraic map $\phi : \oC \rightarrow \TD$ is a \textbf{generalised Harnack curve} in $\TD$ if $\phi$ has degree $\Delta$ and if the logarithmic Gauss map $ \gp : \oC \rightarrow \cp{1}$ is totally real, that is 
\[ \gp^{-1}\big(\rp{1}\big) = \mathbb{R} \oC. \] 
\end{Definition}

By Lemma \ref{lem:critloc}, we have the following reformulation of Definition \ref{defsimpleHarnack}.

\begin{Proposition}\label{gaussmaptotreal}
A real algebraic map $\phi : \oC \rightarrow \TD$ is a generalised Harnack curve if and only if $\phi$ has degree $\Delta$ and the inclusion $\mathbb{R} \oC \subset \oF$ is an equality.
\end{Proposition}

According to Theorem  \ref{thmPR}, Definition \ref{defsimpleHarnack} generalises the Definitions \ref{defmikh1} of smooth simple Harnack curves. It is not hard to see that singular simple Harnack curves of Definition\ref{defmikh2} admit a characterisation similar to Theorem  \ref{thmPR} and that they are particular instances of generalised Harnack curves. 

Observe that there is no restriction on the singularities of the curve $\phi(\oC) \subset \XD$ in Definition \ref{defsimpleHarnack}. We will see later that generalised Harnack curves admit other singularities than real isolated double points.

\begin{Remark}
As a direct consequence of the definition, the real algebraic curve $\oC$ is of type $1$ for any generalised Harnack curve $\phi : \oC \rightarrow \TD$, that is $\oC \setminus \R \oC$ has two connected components.
\end{Remark}

\begin{Theorem}\label{thm:max}
Let $\phi : \oC \rightarrow \TD$ be an immersed generalised Harnack curve. Then, the curve $\oC$ is an M-curve, that is $b_0(\R \oC)=g+1$ where $g$ is the genus of $\oC$.
\end{Theorem}

\begin{proof}
Consider the decomposition $\oC= \C_{>0} \cup \R \oC \cup \C_{<0}$ where $\C_{>0}$ is the locus where $\Ap$ is orientation preserving. Denote $\C_{\geq 0}:=\C_{>0} \cup \R \oC$. Since $\R \oC=\partial \C_{\geq 0}$, the real part $\R \oC$ inherits an orientation from the surface $\C_{\geq 0}$. Denote by $c_0,\dots c_k$ the connected components of $\R \oC$. 
Inside $\C_{\geq 0}:=\C_{>0} \cup \R \oC$, consider a small deformation $\ell$ of $\R \oC$ that coincide with $\R \oC$ outside of an $\varepsilon$-neighbourhood of $\C_\infty$ such that $\ell$ avoids $\C_\infty$ and that the restriction of $\Ap$ to $\ell$ is an immersion. Denote by $\ell_0, \dots, \ell_k$ the corresponding connected components of $\ell$. We carry the orientation of $c_i$ to $\ell_i$. For any immersed oriented curve $c \looparrowright \R^2$, denote the rotational index $\ind(c)$, that is the degree of the map
\[
\begin{array}{rcl}
c & \rightarrow & S^1\\
t & \mapsto & \dfrac{c'(t)}{\vert c'(t) \vert}
\end{array}
\]
where $S^1$ is oriented clockwise. According to \cite[Lemma 2.2]{Kau}, we have that
\begin{equation}\label{eq:kauf}
\ind\big(\Ap(\ell_0)\big)+\dots+\ind\big(\Ap(\ell_k)\big)=k-1.
\end{equation}
Notice also that for any connected component $\ell_j$ corresponding to a component $c_j \subset \R \oC$, we have 
\begin{equation}\label{eq:anycc}
\deg (\gp)_{\vert \ell_j} = 2\cdot \ind\big(\Ap(\ell_j)\big) + \#(\C_{\infty}\cap c_j).
\end{equation}
Indeed, the map $(\gp)_{\vert \R \oC}$ is the composition of $\Ap$ with the above map $t\mapsto\frac{c'(t)}{\vert c'(t) \vert}$ and the $2$-to-$1$ covering $S^1\rightarrow\rp{1}$ away from $\C_{\infty}$. Also, for any point $p\in \C_{\infty}$, there is a unique disc $D_p\subset \C_{\geq 0} \setminus \ell$ adjacent to $p$ and the Gauss map $\Ap(\partial D_p)\rightarrow \rp{1}$ has degree $1$. This justifies the correction term $\#(\C_{\infty}\cap c_j)$ in \eqref{eq:anycc}. From Definition \ref{defsimpleHarnack}, Proposition \ref{degloggauss} and equations \eqref{eq:kauf} and \eqref{eq:anycc}, we deduce that 
\[\deg \gp=-\chi(\C)=2(k-1)+\#\C_\infty \Leftrightarrow k=g.\]
Since the latter is equivalent to $\oC$ being an M-curve, the result follows.
\end{proof}

Let us derive an other characterisation of generalised  Harnack curves from \cite[Appendix A]{MO}. Recall the map 
\[
\begin{array}{rcl}
\Al \, : \, \ttor 		& \rightarrow	& (S^1)^2		\\
(z,w)								& \mapsto		& \big( 2 \, \arg(z), 2 \, \arg(w) \big)
\end{array}
\]
and define $\Alp := \Al \circ \, \phi$ on $\C$. Define $\Cb $ to be the real blow-up of $ \oC $ at every point of $\C_\infty$. Denote by $P_p\subset \Cb$ the fiber of the blow-up over $p \in \C_\infty$. In particular, we have $P_p\simeq\rp{1}$. By construction, we have the 
factorisation $ \C_{\Arg} \rightarrow \Cb \rightarrow \oC$ inducing a double covering $S_p \rightarrow 
P_p$ for any $ p \in \C_{\infty}$. Note that the complex conjugation on $\oC$ extends to an involution on the real surface $\Cb$. The fixed locus of the latter involution consists of the strict transform of $\R \oC$ plus one isolated point in each fibre $P_p$. We denote by $\R \Cb \subset \Cb$ the strict transform of $\R \oC$.

\begin{Lemma} The map $\Alp$ extends to a map $\Alp \, : \, \Cb	\rightarrow (S^1)^2.$
\end{Lemma}

\begin{proof} By Lemma \ref{extarg}, the map $\Alp$ extend to $\C_{\Arg}$. For any $p \in \C_\infty$ and any point $q \in P_p$, the two preimages of $q$ in $S_p$ are mapped to the same value by $\Alp$. Hence, the map $\Alp \, : \, \TC_{\Arg} \rightarrow \sone$ factorises through $\Cb$ and the result follows.
\end{proof}

Define the subset $\Co \subset \Cb$ by  $\Co := \Alp^{-1} \big((1,1)\big)$. Note that $\Al^{-1}\big((1,1)\big)= (\R^*)^2$. Thus, the subset $\Co \subset \Cb$ is the union of 
$\R \Cb$ with some isolated points. The isolated points of $\Co$ come either from
isolated singularities of $\phi(\oC)$ or from non transverse intersection
with a toric divisor. Indeed, for any point $p\in \C_\infty$ such that $\phi(\oC)$ intersects a toric divisor with multiplicity $m$ at $\phi(p)$, there are exactly $m$ point in $P_p\subset \Cb$ belonging to $\Co$ and exactly one of them belongs to $\R \Cb$, see Lemma \ref{extarg}.

Define $ T$ to be the real blow-up of $\sone$ at $(1,1)$ and $ \Cbb $ to be the real blow-up of $ \Cb $ at $ \Co$. As blowing-up at a smooth submanifold of codimension
1 does not change the surface, blowing-up is effective only at isolated points
of $\Co$. The theorem below is the analogue of \cite[Theorem A1]{MO}.

\begin{Theorem}\label{thmalga}
A real algebraic immersion $\phi : \oC \rightarrow \TD$ is a generalised  Harnack curve if and only if $\phi$ has degree $\Delta$ and the map $\Alp \, : \, \Cb \setminus \Co \rightarrow (S^1)^2 \setminus (1,1)$ extends to a covering
\[ \Alp \, : \, \Cbb	 \rightarrow T.\]
\end{Theorem}

\begin{proof}
Assume that $\phim$ is an immersed generalised Harnack curve. By Lemma \ref{lem:critloc} and Proposition \ref{gaussmaptotreal}, the restriction of $\Alp$ to $\C \setminus \R \C$ is a local covering onto its image. As in the proof of \cite[Theorem A1]{MO}, we defined the extension of $\Alp$ at a point $p \in \R \Cb$ to be the image of the real line $\ker T_p \Ap$ by $T_p \Alp$. The image line sits in the tangent space of $\sone$ at $(1,1)$. Equivalently, this line is a point of $T$ in the fibre over $(1,1)\in\sone$. According to the proof Lemma \ref{lem1}, the latter extension is essentially given by the logarithmic Gauss map $\gp$ which is monotonic. Hence, this extension is a covering in a neighbourhood of $\R \Cb$. For any isolated point in $\Co$, the extension of $\Alp$ is given by the blow-up itself, both at the source and at the target. It follows that $\Alp$ extends to a covering $\Alp \, : \, \Cbb	 \rightarrow T$.

Conversely, if $\phi : \oC \rightarrow \TD$ is a real algebraic immersion of degree $\Delta$ such that $\Alp$ extends to a covering $\Alp \, : \, \Cbb	 \rightarrow T$, then $\oF$ is necessarily included in $\R \oC$. According to Proposition \ref{gaussmaptotreal}, the map $\phi$ is a generalised  Harnack curve.
\end{proof}

\begin{Corollary}\label{coralga}
Assume that $\phi : \oC \rightarrow \TD$ is an immersed generalised  Harnack curve. Then
\[ \Areap (\C) = \pi^2 \big(- \chi ( \Cbb ) \big).\]
Assume moreover that $\phi(\oC)$ has no real isolated singularity and cogenus $\delta$. Then 
\[ \Areap (\C) = 2 \pi^2 \big(\Area(\Delta)-\delta). \]
\end{Corollary}
\noindent
\begin{proof}
By Definition \ref{def:areap}, we have that $\int_\C \Alp^{\; \ast} (dx_2\wedge dy_2)= 4 \cdot \Areap (\C) $. 
By Theorem \ref{thmalga}, we have that $\Alp$ extends to a covering $\Alp \, : \, \Cbb	 \rightarrow T$ of degree $-\chi( \Cbb)\big/\chi( T)=-\chi( \Cbb)$. As the blown-up torus $T$ coincide with $\sone$ outside of a set of measure $0$, we have
$$ \int_\C \Alp^{\; \ast} (dx_2\wedge dy_2)= \deg(\Alp) \cdot \int_{\sone}dx_2\wedge dy_2 = 4\pi^2\cdot\big(-\chi( \Cbb)\big).$$ 
The first part of the statement is proven. For the second part, we need to prove that $\chi ( \Cbb )=2\delta-\Area(\Delta)$. Under the present assumptions, the surface $\Cbb$ is the real blow-up of $\oC$ at $b_\Delta$-many points. It follows that 
\[ \chi(\Cbb)=\chi(\oC)-b_\Delta= \big(2-2(g_\Delta-\delta)\big)-b_\Delta= 2\delta-\big(2g_\Delta+b_\Delta-2\big)=2\big(\delta-\Area(\Delta)\big) \]
by Pick's Formula.
\end{proof}

\section{Tropical constructions}\label{sec:trop}

\subsection{Phase-tropical curves}

\tocless\subsubsection{Tropical curves}

Let us recall some standard notions about tropical curves in the plane. All definitions, statements and their proofs can be found in \cite{Mikh05}, \cite{IMS},  and \cite{BIMS}.\\
A \textbf{tropical Laurent polynomial} in two variables $x$ and $y$ is a function
\[f(x,y) =`` \sum_{(\alpha, \beta) \in A} c_{(\alpha, \beta)} x^\alpha y^\beta"\]
where $A \subset \mathbb{Z}^2$ is a finite set and the usual arithmetic operations are replaced by the tropical ones
\[ ``x+y":= \max \left\lbrace x,y \right\rbrace \, \text{ and } \,  ``xy":= x+y. \]
Any tropical Laurent polynomial is piecewise affine linear and convex. The \textbf{Newton polygon} $\New(f)$ of $f$ is the convex hull of $A$ in $\mathbb{R}^2 = \mathbb{Z}^2 \otimes_\mathbb{Z} \mathbb{R}$. The \textbf{tropical zero set} $\cZ(f)$ of $f$ is defined as the subset of $\mathbb{R}^2$ where $f$ is not smooth. Formally, we have $$\cZ(f)=\left\lbrace(x,y) \in \R^2 \, \big| \, f(x,y)=``c_{(\alpha_1, \beta_1)} x^{\alpha_1} y^{\beta_1}"=``c_{(\alpha_2, \beta_2)} x^{\alpha_2} y^{\beta_2}", \; (\alpha_1, \beta_1)\neq (\alpha_2, \beta_2) \in A \right\rbrace.$$
Consequently, any non-empty tropical zero set  is a piecewise linear graph with rational slopes in $\mathbb{R}^2$.
If $g$ is another tropical Laurent polynomial given by $$g(x,y) := ``c_{(\alpha, \beta)} x^\alpha y^\beta \cdot f(x,y)",$$ then $\cZ(f)=\cZ(g)$ but the converse fails to be true.
For a tropical Laurent polynomial $f$ with Newton polygon $\Delta$, consider its \textbf{extended Newton polygon} 
\[ \tilde{\Delta} : = \conv \left\lbrace \big( (\alpha, \beta) , t \big) \in \mathbb{R}^3  \: \big| \:  (\alpha, \beta) \in A, \: t \geq c_{(\alpha, \beta)} \right\rbrace. \]
The projection on the first two coordinates of the union of all closed bounded faces of $\tilde{\Delta}$ induces a subdivision $\Sub_f$ of $\Delta$. The Legendre transform gives rise to the following duality, see \cite[Theorem 3.3]{IMS}.

\begin{Proposition}
Let $f$ be a tropical Laurent polynomial in two variables. The subdivision of \,$\mathbb{R}^2$ induced by $\cZ(f)$ is dual to the subdivision $\Sub_f$ of $\Delta$ in the following sense

-- the 2-cells of $\, \mathbb{R}^2 \setminus \cZ(f)$ are in bijection with vertices of $\, \Sub_f$ and the 2-cells of $\, \Sub_f$ are in bijection with the vertices of $\, \cZ(f)$,

-- any edge of $\, \cZ(f)$ is in bijection with an edge of $\, \Sub_f$ with  orthogonal direction.

-- The above correspondence reverses the incidence relation.

Moreover, the unbounded 2-cells of $\, \mathbb{R}^2 \setminus \cZ(f)$ are dual to boundary points of $\Delta$ and unbounded edges of $\, \cZ(f)$ are dual to edges on the boundary of $\Delta$.
\end{Proposition}

\begin{Definition}\label{deftropcurv}
Let $f$ be a tropical Laurent polynomial in two variables. For any edge $\varepsilon$ of $\cZ(f)$, the \textbf{weight} $w(\varepsilon)$ of $\, \varepsilon$ is the integer length of its dual edge $\varepsilon^\vee$ in $\Sub_f$, that is $\#( \varepsilon^\vee \cap \mathbb{Z}^2) -1$.
A \textbf{tropical curve} $C \subset \mathbb{R}^2$ is a tropical zero set $\cZ(f)$ together with the weights $w(\varepsilon)$.
If $\Delta$ is the Newton polygon of $f$, denote by $\Sub_C:=\Sub_f$ the \textbf{subdivision} of $\Delta$ dual to $C$.
\end{Definition}

\begin{Remark}
The convex piecewise affine linear function $S_f$ defining the spine of the curve $\left\lbrace f=0 \right\rbrace \subset \ttor$ is a tropical Laurent polynomial, see Section \ref{sec:loggeo}. Equipped with the collection of weights of Definition \ref{deftropcurv}, the spine of an algebraic curve in $\ttor$ is a tropical curve.
\end{Remark}

\begin{Definition}\label{stcnorm}
An \textbf{abstract tropical curve} is a finite graph with $1$-valent vertices removed and equipped with a complete inner metric. A proper continuous map $h:C\rightarrow \R^2$ from an abstract tropical curve $C$ is a \textbf{parametrised tropical curve} if the image of any unit tangent vector to $C$ under the differential $d h$ is in $\Z^2\subset\R^2$, and if for each vertex $v\in C$ we have 
\[\sum_e u(e)=0\]
where the sum is taken over all edges adjacent to $v$, and $u(e)$ is the image under $d h$ of the unit tangent vector to $e$ such that this vector points outward of $v$.\\
A parametrised tropical curve $h:C\rightarrow \R^2$ is \textbf{nice} if the image of any unit tangent vector to $C$ under the differential $d h$ is a primitive integer vector and if $h$ is injective outside of a finite subset of $C$. The parametrised tropical curve $h:C\rightarrow \R^2$ is \textbf{very nice} (\VS for short) if moreover $C$ is trivalent.
\end{Definition}

\begin{Remark}
The image $h(C)\subset \R^2$ of a parametrised tropical curve is a tropical curve in $\R^2$ when equipped with the appropriate weights. For nice parametrised tropical curves, these weights are all equal to $1$.
\end{Remark}

\begin{Definition}\label{def:tropedgeleafnode}
Let $h:C\rightarrow \R^2$ be a parametrised tropical curve. We denote the set of vertices and edges of $C$ by $V(C)$ and $E(C)$ respectively. As a convention, edges are always open. If $h$ is nice, the points of $h(C)$ having several preimages in $C$ are called the \textbf{nodes} of $C$. The \textbf{multiplicity} of a node $n$ is the integer $m(n) := 2\Area (n^\vee)$
where $n^\vee$ is the 2-cell dual to $n$ in $\Sub_{h(C)}$. A node $n$ is \textbf{hyperbolic} if $m(n)=2$. In particular, it has exactly two preimages in $C$.\\
The \textbf{Newton polygon} of $h:C\rightarrow \R^2$ is the Newton polygon of any tropical Laurent polynomial defining $h(C)$. It is defined only up to translation.
\end{Definition}

We end up this subsection by recalling what is the stable intersection of two tropical curves $C_1$ and $C_2$ in $\R^2$. The stable intersection of $C_1$ and $C_2$ is a formal sum of point in $C_1 \cap C_2$. If $C_1$ and $C_2$ intersect each other transversally (in particular, away from their vertices), then the stable intersection of $C_1$ and $C_2$ is the sum
\[\sum_{p \in C_1 \cap C_2} m_p(C_1,C_2) \cdot p \] 
where $m_p(C_1,C_2)= 2 \Area p^\vee$ and $p^\vee$ is the 2-cell dual to $p$ in $\Sub_{C_1 \cup C_2}$.
For any curves $C_1$ and $C_2$, there is an open dense subset 
$\mathscr{U} \subset \R^2$ such that for any $\vec{v} \in \mathscr{U} $, the tropical curves $C_1$ and $\C_2 + \vec{v}$ intersect transversally as above. The stable intersection of $C_1$ and $C_2$ is the limit of the stable intersection of $C_1$ and $C_2+ \vec{v}$ when $\vec{v}$ tends to $0$ and it does not depend on $\vec{v}$,  see \cite{RGST} for more details.

\tocless\subsubsection{Phase-tropical curves}\label{sec:ptcurve}

\begin{figure}[h]
\centering
\scalebox{1}{\input{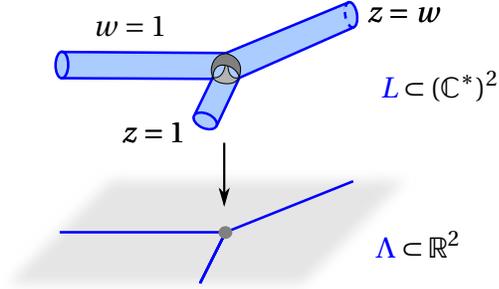}}
\caption{The fibration $ \A : L \rightarrow \Lambda$.}
\label{fig:L0}
\end{figure}

\begin{figure}[h]
\centering
\input{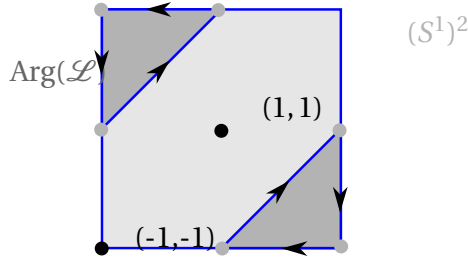}
\caption{The coamoeba of $\mathcal{L}$ (in dark gray), and the framing of its 3 boundary geodesics (in blue).}
\label{figcoamoeba}
\end{figure}

As in \cite[\S 6]{Mikh05}, consider the change of the holomorphic structure
\[
\begin{array}{rcl}
H_t \: : \: \ttor & \rightarrow & \ttor \\
(z,w) & \mapsto & \displaystyle \left( \vert z \vert^{\frac{1}{\log(t)}} \frac{z}{\vert z \vert} , \vert w \vert^{\frac{1}{\log(t)}} \frac{w}{\vert w \vert} \right)
\end{array}.
\]
Observe that 
\[\mathcal{A} \circ H_t = \frac{1}{\log(t)} \mathcal{A}.\]
Denote $ \mathcal{L} := \left\lbrace (z,w) \in \ttor \: \big| \: z+w+1=0  \right\rbrace$. The sequence of topological surfaces $\big\lbrace H_t(\mathcal{L}) \big\rbrace_{t >1}$ converges in Hausdorff distance to the so-called phase-tropical line $L$ when $t$ tends to $\infty$. Let us describe  $L$. The amoeba $\A(L)\subset \R^2$ is the tropical line $\Lambda:=\cZ(``0+x+y")$ consisting of $3$ half-rays emanating from $(0,0)$ and directed by $(-1,0)$, $(0,-1)$ and $(1,1)$ respectively. The preimage by $\A_{\vert L}$ of the latter open rays are  cylinders in $\{z=1\}$, $\{w=1\}$ and $\{z=w\}$ respectively. These cylinders glue to the preimage by $\A_{\vert L}$ of $(0,0)$: the latter is the pair of pants given  as the closure of $\Arg(\mathcal{L}) \subset \sone$. The phase-tropical line $L$ is then homeomorphic to  a sphere with $3$ punctures, see Figures \ref{fig:L0} and \ref{figcoamoeba}.

Here, a \textbf{toric transformation} is a map  $ A : \ttor \rightarrow \ttor $ of the form 
\[ (z,w) \mapsto \big( z_0 z^{a}w^{b}, w_0 z^{c}w^{d} \big) \]
where $(z_0,w_0) \in  \big( \mathbb{C}^\ast \big)^2 $ and $  \begin{psmallmatrix} a & b \\ c & d \end{psmallmatrix} \in \text{GL}_{2} \big( \mathbb{Z} \big)$.
The map $A$ descends to an affine linear transformation on $\mathbb{R}^2$ (respectively on $\sone$) by composition with the projection $\A$ (respectively $\Arg$) that we still denote by $A$.

\begin{Definition}
A \textbf{general phase-tropical line} in $\ttor$ is the image of $L$ by any toric transformation. 
\end{Definition}

For a $3$-valent abstract tropical curve $C$ and any vertex $v \in V(C)$, we denote by $Y_v \subset C$ the \textbf{tripod} obtained as the union of $v$ with its $3$ adjacent edges in $C$.

\begin{Definition}\label{sctcdef} 
A \textbf{simple phase-tropical curve} $V \subset \ttor$ supported on a \VS parametrised tropical curve $h:C\rightarrow \R^2$ (see Definition \ref{stcnorm}) is a topological space such that :

-- $ \mathcal{A} (V)=h(C) $,

-- for any $v \in V(C)$, there exists a unique general phase-tropical line $\Gamma_v \subset \ttor$ such that $h(Y_v) \subset \A(\Gamma_v)$ and $\Gamma_v \cap \A^{-1}\big(h(Y_v)\big) \subset V$,

-- for any $e \in E(C)$, $v_1$ and $v_2$ its two adjacent vertices in $C$, $\Gamma_{v_1}$ and $\Gamma_{v_2}$ coincide on $\mathcal{A}^{-1} \big((h(e)\big)$,

-- $V= \bigcup_{v\in V(C)}  \Gamma_v \cap \A^{-1}\big(h(Y_v)\big)$.
\end{Definition}

\begin{Remark}\label{rem:normpt}
By construction, a simple phase-tropical curve $V \subset \ttor$ is a smooth topological surface away from finitely many transversal self-intersection points sitting over the nodes of tropical curve $h(C) \subset \R^2$. The latter surface can be normalised by blowing up $\ttor$. The result is a smooth topological surface of genus $g$ with $n$ puncture where $g:=b_1(C)$ and $n$ is the number of unbounded edges of $C$. This normalisation can be obtained as a phase-tropical morphism but we do not need this formalism here. We refer to \cite[\S 6.3]{Mikh05} and \cite{L}.
\end{Remark}

Similarly to Riemann surfaces, simple phase-tropical curves can be described in terms of Fenchel-Nielsen coordinates, see \cite{L}. Recall that $ \Arg(\mathcal{L})$ is the union of the $2$ open triangles delimited by the $3$ geodesics $\{\arg(z)=-1\}$, $\{\arg(w)=-1\}$ and $\{\arg(z)=-\arg(w)\}$ (in multiplicative notation on $\sone$) plus their 
$3$ common vertices. On each of the boundary geodesics of $\Arg (\mathcal{L})$, we fix an orientation as shown on Figure \ref{figcoamoeba}. We fix the origin of the geodesics $\{\arg(z)=-1\}$, $\{\arg(w)=-1\}$ and $\{\arg(z)=-\arg(w)\}$ to be $(-1,-1)$, $(1,-1)$ and $(-1,1)$ respectively. The orientation of each geodesic as well as the choice of the origin is preserved by each of the $6$ toric transformations preserving $\mathcal{L}$. Each pointed and oriented geodesic is canonically isomorphic to $S^1$ as an Abelian group.

For a general phase-tropical line $\Gamma := A(L)$, we have that 
$\A (\Gamma) = A (\Lambda) $ and that the fiber in $\Gamma$ over the vertex of $\A (\Gamma)$ is $\Arg (\Gamma) = A \big( \Arg(L) \big) = A \big( \overline{\Arg(\mathcal{L})} \big) $.
We carry the orientation and origin of the boundary geodesics of $\Arg(L)$ to the boundary geodesics of $\Arg (\Gamma) = A \big( \Arg(L) \big)$ using the map $A$. This construction does not depend on the choice of $A$. With the latter structure, the boundary geodesics of $\Arg (\Gamma)$ are canonically isomorphic to $S^1$ as well.
 
For a simple phase-tropical curve $V$ with underlying \VS parametrised tropical curve $h:C\rightarrow \R^2$, and two vertices $v_1, v_2 \in V(C)$ connected by an edge $e \in E(C)$, the holomorphic annulus $\left( \mathcal{A}_{|_V} \right)^{-1} \big(h(e)\big) $ maps to a geodesic $\gamma_e$ in the argument torus. This is a common boundary geodesic of the two coamoebas $ \left( \mathcal{A}_{|_V} \right)^{-1} \big((h(v_1)\big) $ and 
$ \left( \mathcal{A}_{|_V} \right)^{-1} \big((h(v_2)\big) $. Hence, the geodesic $\gamma_e$ inherits two identifications $\tau_1, \tau_2: S^1 \rightarrow \gamma_e$ from the above construction such that $\tau_2^{-1} \circ \tau_1$ is of the form $z  \mapsto  -\overline{e^{i\theta}z}.$ Since $\tau_2^{-1} \circ \tau_1$ is involutive, the element $e^{i\theta} \in S^1$ does not depend on the ordering of $v_1$ and $v_2$.

\begin{Definition}\label{twistdef}
Let $V\subset \ttor$ be a simple phase-tropical curve supported on $h:C\rightarrow \R^2$. For any bounded edge $e \in E(C)$, the \textbf{twist parameter} of the edge $e$ is the element $e^{i\theta} \in S^1$ constructed above.
\end{Definition}

\begin{Definition}
A \textbf{simple real-tropical curve} $ V \subset \ttor $ is a simple phase-tropical curve that is globally invariant under complex conjugation. We denote by $ \mathbb{R} V \subset ( \mathbb{R}^\ast)^2 $ its real point set and by $ \mathbb{T} V $ the image of \,$ \mathbb{R} V $ under the diffeomorphism
\[\begin{array}{rcl}
\mathcal{A}_s \, : \, ( \mathbb{R}^\ast)^2 & \rightarrow & \mathbb{R}^2 \times \{\pm 1\}^2 \\
(x,y) & \mapsto & \Big( \big(\frac{x}{\vert x \vert}\ln \vert x \vert , \frac{y}{\vert y \vert} \ln \vert y \vert\big), \big(\frac{x}{\vert x \vert},\frac{y}{\vert y \vert}\big) \Big).
\end{array}\]
\end{Definition}

\begin{Remark}
Observe that the amoeba map $\A$ is equal to the composition of the map $ \mathcal{A}_s $ with 
\[\begin{array}{rcl}
\Abs \, : \, \mathbb{R}^2 \times \{\pm 1\}^2  & \rightarrow &  \mathbb{R}^2 \\
\big( (x,y),(\sigma,\nu) \big) & \mapsto & (\sigma x, \nu y),
\end{array}\]
Observe  moreover that for any $e \in E(C)$, the locus of the holomorphic annulus $ \left( \mathcal{A}_{|_V} \right)^{-1} \big((h(e)\big) $ fixed by complex conjugation has  two connected components when the annulus is real and is empty otherwise.
\end{Remark}

The following results are easy consequences of the above definitions. The proofs are left to the reader.

\begin{Proposition}\label{prop2to1}
For a simple real tropical curve $V \subset \ttor $ supported on $h:C\rightarrow\R^2$, the set $\T V$ is a piecewise linear curve. Moreover, the restriction $\Abs \, : \, \T V \rightarrow C$ is 2-to-1 and maps any maximal domain of linearity of $\T V$ onto $h(e)$ for some $e\in E(C)$.
\end{Proposition}

\begin{Proposition}
The twist parameters of a real-tropical curve $ V \subset \ttor $ are always in $ \left\lbrace -1, 1 \right\rbrace \subset S^1$.
\end{Proposition} 

\begin{Definition}\label{def:leafedgept}
For a simple real-tropical curve $V \subset \ttor$ supported on $h:C\rightarrow\R^2$, an \textbf{edge} of \,$ \T V $ is any of the two connected components of the preimage by $\A_s$ of $h(e)$ for some $e\in E(C)$. We denote the set of edges of \,$ \T V $ by $E(\T V)$. An \textbf{inflection pattern} of \,$ \mathbb{T} V $ is a collection of three consecutive edges $e_1$, $e_2$ and $e_3$ of \,$\mathbb{T} V$ such that the piecewise linear arc $e_1 \cup e_2 \cup e_3$ is not convex. Finally, a bounded edge $e \in E(C)$ is said to be \textbf{twisted} if its twist parameter is -1.
\end{Definition}

\begin{Lemma}\label{leminfpat}
Let $V \subset \ttor$ be simple real-tropical curve supported on $h:C\rightarrow\R^2$. An edge $e\in E(\T V)$ is the middle edge of inflection pattern if and only if the edge $e'\in E(C)$ such that $\Abs(e)=h(e')$ is twisted. In particular, the inflection patterns of $ \mathbb{T} V $ are in 2-to-1 correspondence with the twisted edges of $C$.
\end{Lemma}

\begin{figure}[h]
\centering
\input{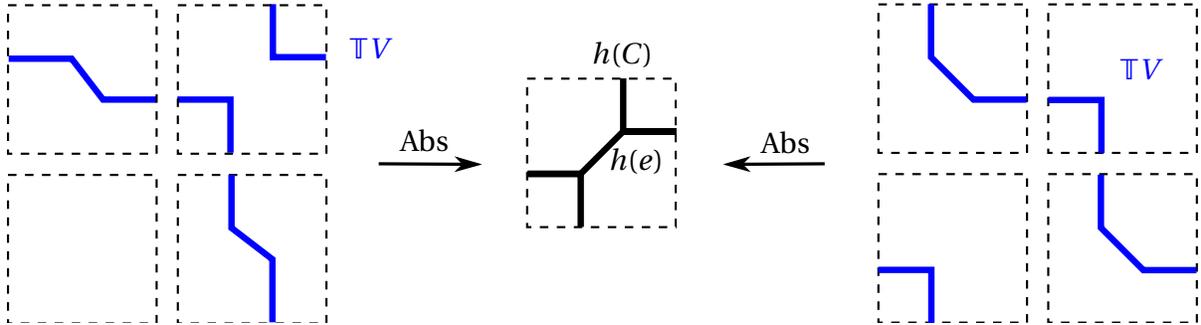}
\caption{The map $\Abs:\T V\rightarrow h(C)$ in a neighbourhood of $h(e)$, $e\in E(C)$, in the twisted case (left) and in the non-twisted case (right). Observe the two inflections patterns mapping to $h(e)$ in the twisted case.}
\label{twist}
\end{figure}

\begin{proof}
The latter can be formulated in terms of signs distribution in combinatorial patchworking, see for example section 3 of  \cite{BIMS} and references therein. A simple computation allows to describe the two possible cases pictured in Figure \ref{twist}.
\end{proof}

\subsection{Tropical  Harnack curves}\label{subsecthc}

Let $h:C \rightarrow \mathbb{R}^2 $ be a nice parametrised tropical curve. Consider a loop $\lambda \subset C$. Note that $\lambda$ can be considered as a subset of $E(C)$. Denote by $ \Gamma_\lambda \subset E(C)$ the subset of edges $e$ in $\lambda$ such that the piecewise linear arc in $\R^2$ formed by $h(e)$ and the image of its two adjacent edges in $\lambda$ is not convex.

\begin{Definition}\label{deftropsimpleHarnack}
A \textbf{tropical Harnack curve} $h: C \rightarrow \mathbb{R}^2 $ is a nice parametrised tropical curve  such that for every loop $\lambda \subset C$, one has 
\begin{equation}\label{thc}
\sum_{e \in \Gamma_\lambda} u(e)\equiv 0  \mod  2
\end{equation}
where $u(e)$ is the image by $d h$ of any unit tangent vector on $e$.
\end{Definition}

\begin{figure}[h]
\centering
\input{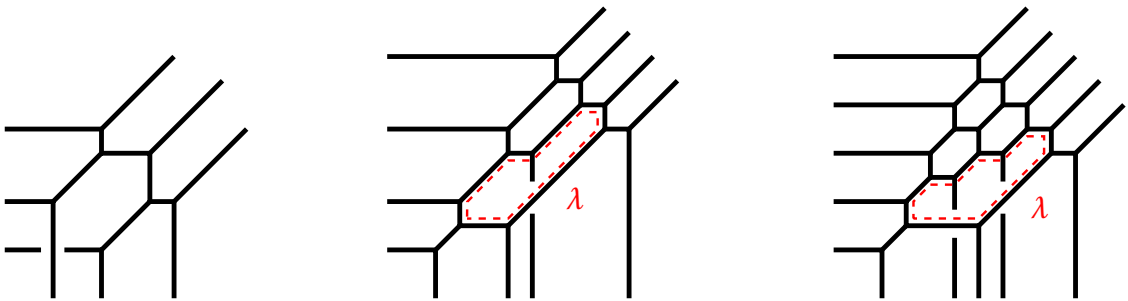}
\caption{Three \VS parametrised tropical curves: a cubic, a quadric and quintic from left to right. The cubic and the quintic are tropical Harnack curves while the quartic is not. The cubic is rational so that the condition of Definition \ref{deftropsimpleHarnack} is empty. The red loop $\lambda$ in both the quartic and the quintic is the only one for which $\Gamma_\lambda$ is non-empty.}
\label{TropCurves1}
\end{figure}

\begin{Remark}
For any loop $\lambda$ as in Definition \ref{deftropsimpleHarnack}, the set $\Gamma_\lambda$ is non-empty if and only if $h(\lambda)$ contains a node of $h(C)$, see Definition \ref{def:tropedgeleafnode} and Figure \ref{TropCurves1} for examples.
\end{Remark}

\begin{Proposition}\label{equivdeftropsimpleHarnack}
A \VS parametrised tropical curve $h:C\rightarrow \mathbb{R}^2 $ is a tropical Harnack curve if and only if there exists a simple real-tropical curve $V \subset \ttor$ supported on $h$ and such that such that $\T V$ has no inflection pattern.
Moreover, the phase-tropical curve $V$ is unique up to the four sign changes of the coordinates in $\ttor$. 
\end{Proposition}

\begin{Definition}\label{def:ptharnack}
A \textbf{phase-tropical Harnack curve} $V \subset \ttor$ is a simple real-tropical curve supported on a tropical Harnack curve such that $\T V$ has no inflection pattern.
\end{Definition}

\begin{figure}[h]
\centering
\input{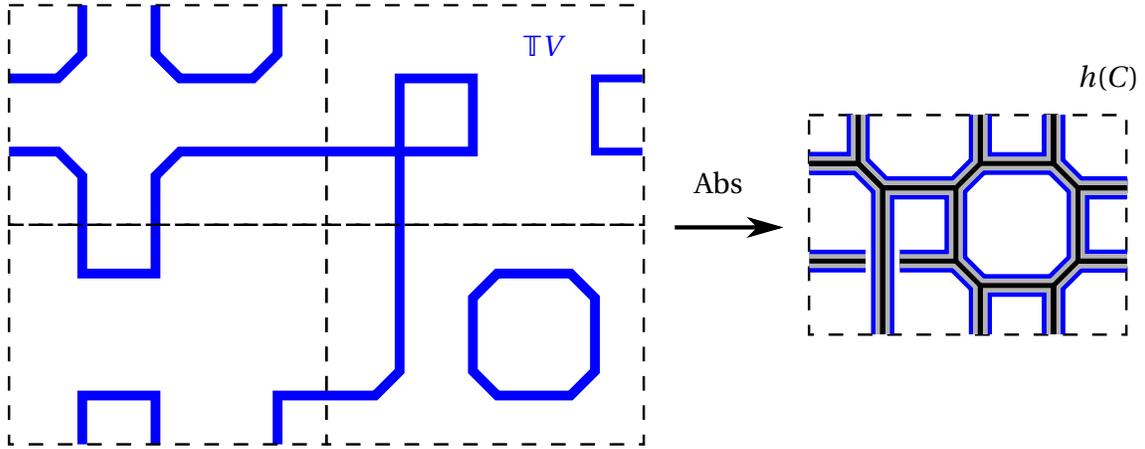}
\caption{Recovering $\T V$ from the topical Harnack curve $h:C\rightarrow\R^2$. On the right, the immersed ribbon $R$ is pictured in gray with boundary in blue.}
\label{TropCurves3}
\end{figure}

Before proving Proposition \ref{equivdeftropsimpleHarnack},  let us recall how one can recover the curve $\T V$ for a phase-tropical Harnack curve $V$ from its underlying tropical  Harnack curve $ h:C \rightarrow\R^2$. According to Propositions \ref{prop2to1} and \ref{equivdeftropsimpleHarnack}, each connected component of $\T V$ is convex and so is its image by $\Abs$. Also, there are exactly two connected components of $\T V$ whose image by $\Abs$ contains a given element in $E(\T V)$. Now, consider an infinitely thin ribbon $R$ containing  $C$ and such that the parametrisation $h: C\rightarrow \R^2$ extends to an immersion on $R$. By construction, the image by  $h$ of any connected component of $\partial R$ is convex and to any edge in $E(\T V)$ correspond exactly two connected components of $\partial R$. Up to sign, the curve $\T V $ can be recovered as follows. Each connected component $c_V \subset \T V$ corresponds to a connected component $c_R \subset \partial R$ in such a way that  $h(c_R)=\Abs(c_V)$. For a given edge $e \in E(C)$, the two components of $\T V$ whose image contain $e$ sit in two different quadrants of $\mathbb{R}^2 \times \{\pm 1\}^2$ given by signs $(\sigma_1,\nu_1)$ and $(\sigma_2,\nu_2)$ satisfying $(\sigma_1,\nu_1)-(\sigma_2,\nu_2)\equiv (a,b) \mod 2$ where $(a,b)$ is a primitive integer vector supporting $h(e)$. Indeed, the fiber of $\A_{\vert V}$ over $h(e)$ is a holomorphic annulus given by an equation of the form
\[ z^{-b}w^{a}=c\]
where $c \in \mathbb{R}^\ast$. The two connected components of the real part of the latter annulus sit in two different quadrants of $(\mathbb{R}^*)^2 $ whose respective signs $(\sigma_1,\nu_1)$ and $(\sigma_2,\nu_2)$ satisfy $(\sigma_1,\nu_1)-(\sigma_2,\nu_2)\equiv (a,b) \mod 2$. All the latter information allow then to recover $\T V $ up to sign, since one needs to know how to lift at least one of the components of $\partial R$. This construction is illustrated in Figure \ref{TropCurves3}.

\begin{Remark}\label{rem:ptHmax}
The normalisation $\tilde V$ of a phase-tropical Harnack curve $V$ inherits a real structure. The real locus $\R \tilde V \subset \tilde V$ is maximal in the sense that $b_0(\R \tilde V)=g+n$ where $\tilde V$ is a topological Riemann surface of genus $g$ with $n$ punctures. 
\end{Remark}

\begin{proof}[Proof of Proposition \ref{equivdeftropsimpleHarnack}.]
Suppose there is a simple real-tropical  curve $V\subset \ttor$ supported on $h:C\rightarrow\R^2$  and such that  $\T V$ has no inflection pattern. By  Lemma \ref{leminfpat}, it is equivalent to the fact that $V$ has no twisted edges.
For any vertex $v \in \lambda$, there is a distinguished point among the three vertices of the coamoeba $ \left( \mathcal{A}_{|_V} \right)^{-1} \big(h(v)\big) $, namely the intersection point of the two geodesics corresponding to the two edges in $\lambda$ adjacent to $v$. Let us look at the position of this distinguished point in the argument torus $(S^1)^2$ while going around $\lambda$.
Going from a vertex $v$ to the next one via an edge $e$, the point is moved according to the following rule : if $e$ is not in $ \Gamma_\lambda$, then this point is fixed; if $e$ is in $ \Gamma_\lambda$, this point is moved by $\pi\cdot u(e)$ in $\sone\simeq (\R/2\pi\Z)^2$. After a full cycle, the distinguished point has to come back to its initial place. This is equivalent to the condition stated in definition \ref{deftropsimpleHarnack} on the loop $\lambda$. Hence $C$ is a tropical  Harnack curve.

Conversely, if $C$ satisfies the condition of definition \ref{deftropsimpleHarnack} for any cycle, pick an initial vertex $v_0$ on $C$. There is exactly one possible general phase-tropical line $\Gamma_{v_0}$ such that $ \left( \mathcal{A}_{|_V} \right)^{-1} \big(h(Y_{v_0})\big) = \Gamma_{v_0} \cap \A^{-1} \big(h(Y_{v_0})\big) $, up to the four changes of signs of the coordinates. 
The twists determine the gluing of the general phase-tropical lines above adjacent vertices along the common edge. Hence, 
once $\Gamma_{v_0}$ is fixed, the adjacent general phase-tropical lines are also fixed. The first part of the proof shows that the condition of definition \ref{deftropsimpleHarnack} is necessary and sufficient for this construction to close up along every cycle $\lambda$. The proposition is proved.
\end{proof}

\subsection{Construction by tropical approximation}

Recall from Remark \ref{rem:normpt} that any simple phase-tropical curve $V \subset \ttor$ supported on  $h:C \rightarrow \R^2$ can be normalised into a smooth topological surface whose topology is governed by $C$.

\begin{theorem}[Mikhalkin]\label{thmmikh}
Let $V \subset \ttor$ be a simple real-tropical curve  supported on $h:C\rightarrow \R^2$ with Newton polygon $\Delta$ such that the normalisation $\tilde{V} $ of $V$ has genus $g$ and $n$ punctures. Then, there exists a family of Riemann surfaces $ \left\lbrace \C_t \right\rbrace_{t>>1} \subset \mathcal{M}_{g,n} $ together with immersions $\phi_t \: : \: \C_t \rightarrow \ttor$
such that 
\begin{enumerate}
\item[$\ast$] $\phi_t \big( \C_t \big)$ is a real algebraic curve with newton polygon $\Delta$,
\item[$\ast$] $\phi_t \big( \C_t \big)$ converges in Hausdorff distance to $V$.
\end{enumerate}
\end{theorem}

\begin{proof}
This theorem is disseminated in \cite{Mikh05}. We refer to  \cite[Theorem 5]{L} for a proof.
\end{proof}

\begin{Definition}
Let $h:C\rightarrow\R^2$ be a tropical Harnack curve with Newton polygon $\Delta$. Define $\Top(h)$ to be the topological type
\[\Big( \R \TD, \; \overline{\R V} \; \bigcup_{j\in \Z/n\Z} \R \mathcal{D}_j \Big)\]
where $V$ is any phase-tropical Harnack curve sitting above $C$ (see Proposition \ref{equivdeftropsimpleHarnack} and Definition \ref{def:ptharnack}).
\end{Definition}

The main result of this section, that will arise as a consequence of Theorem \ref{thmmikh}, is the following.

\begin{Theorem}\label{thmexistharncurv}
Let $h:C\rightarrow\R^2$ be a tropical  Harnack curve with Newton polygon $\Delta$. Then, there exists a generalised Harnack curve $ \C \subset \TD$ such that 
\[\Big( \R \TD, \; \R \C  \bigcup_{j\in \Z/n\Z} \R \mathcal{D}_j \Big) = \Top(h).\]
\end{Theorem}

Before giving the proof, let us illustrate Theorem \ref{thmexistharncurv} with some examples. Recall that the simple Harnack curves considered previously in the literature could only develop real isolated double points as singularities, see \cite{MR}. Using the latter theorem, we are able to construct generalised Harnack curves with different types of singularity.

\begin{figure}[h]
\centering
\input{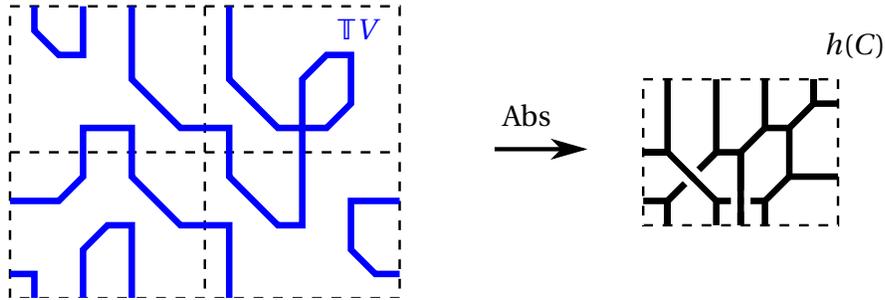}
\caption{$\T V \subset \mathbb{P}^1 \times \mathbb{P}^1$ for the tropical Harnack curve $C$.}
\label{CompConj}
\end{figure}

Consider first the tropical Harnack curve of Figure \ref{TropCurves3}. According to the latter theorem, there exists a generalised Harnack curve $\C \subset \mathbb{P}^1\times \mathbb{P}^1$ whose topological type is described by the curve $\T V$, on the left of Figure \ref{TropCurves3}. In particular, the real part of $\C$ contains an hyperbolic node.  As a second example, consider the tropical Harnack curve of Figure \ref{CompConj}. The algebraic curve $\C \subset \mathbb{P}^1\times \mathbb{P}^1$ provided by Theorem \ref{thmexistharncurv} is  a curve of bi-degree $(4,2)$ in $\mathbb{P}^1 \times \mathbb{P}^1$ with one hyperbolic node and two complex conjugated nodes coming from the two edges of $C$ intersecting each other with multiplicity $4$. 

\begin{figure}[h]
\centering
\input{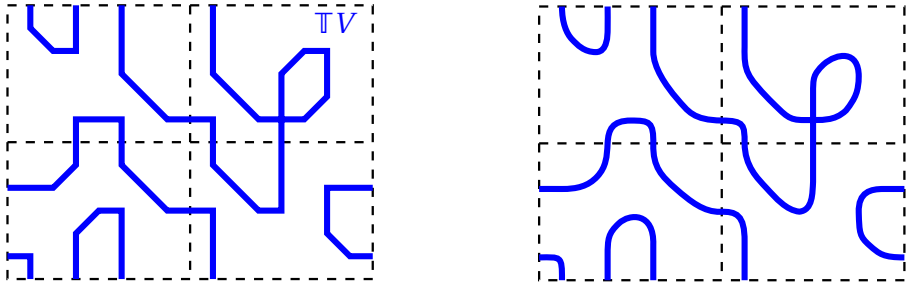}
\caption{A deformation of $\T V$  into a smoothly immersed curve.}
\label{fig:smooth}
\end{figure}

We now proceed to the proof of Theorem \ref{thmexistharncurv}. Let $V$ be a phase-tropical Harnack curve supported on $h:C\rightarrow\R^2$ and with Newton polygon $\Delta$. Our aim is to define an analogue of the real  logarithmic curvature on $\R V$ (see \cite[\S 5]{BLR}). The normalisation $\pi : \tilde V \rightarrow V$ induces a parametrisation of $\R V$ (and then $\T V$) by $\R \tilde V$. We can compactify the latter parametrisation by filling the punctures of $\tilde V$ at the source and considering the closure of $\R V$ in $\XD$. Moreover, we can orient the source according to one of the two halves of $\tilde V \setminus \R \tilde V$. Now, consider any isotopy of the latter parametrisation in $\R\XD$ inducing a smoothing of the vertices of $\T V$, see Figure \ref{fig:smooth}. Define the \textbf{total logarithmic curvature} $\kappa_{\R V}$ of $\R V$ to be the absolute value of the degree of the Gauss map from the deformation of the parametrised curve $\overline{\T V}$ into $\rp{1}$. Clearly, the latter number depends neither on the deformation nor on the choice of orientation on $\R \tilde V$.

\begin{Proposition}\label{proptotcurv}
The total logarithmic curvature $\, \kappa_{\R V}$ is equal to 
$-\chi(\tilde{V})$.
\end{Proposition}

\begin{proof}
On the one hand, notice that $-\chi(\tilde{V})$ is equal to the number of vertices of $C$. Indeed, the surface $\tilde V$ is obtained by gluing pairs of pants, one for each vertex in $ C$.
On the other hand, we can compute the logarithmic curvature of  $\R V$ by computing contribution over every tripod in $C$. Formally, we can choose an $\varepsilon$-deformation of $\T V$ tangent to $\T V$ at the middle of every edge. The middle points of all edges of $\T V$ cut the deformation of $\T V$ into $3\cdot \#V(C)$ arcs. Since $\T V$ has no inflection pattern, we can choose one of the two orientations on $\R \tilde V$ such that all the $3\cdot \#V(C)$ oriented arcs are convex. In particular, they contribute positively to the Gauss map. Now, for any vertex $v$ of $C$, there are exactly $3$ of the latter arcs whose image by $\Abs$ intersect the $\varepsilon$-neighbourhood of $h(v)$. It is easy to see that the contribution to the Gauss map of these $3$ arcs is $1$ (see also \cite[Proposition 5.8]{BLR}). It follows that $\kappa_{\R V}$ is equal to the number of vertices of $C$.
\end{proof}

\begin{proof}[Proof of Theorem \ref{thmexistharncurv}.]
For $t$ large enough, the real part of the immersed curve $\phi_t \big( \C_t \big)$ in Theorem \ref{thmmikh} realises the topological type $\Top(h)$ and the parametrisation  $(\phi_t)_{\vert \R \C_t}$ is a small deformation of the parametrisation of $\overline{\R V}\subset \XD$ considered above. On the one hand, we know that $\kappa_{\R V}$ is given by $-\chi(\tilde{V})$ by the previous proposition. On the other hand, the total logarithmic curvature $\kappa_{\R V}$ is the total logarithmic curvature of $ \phi_t(\R \C_t)$, which is, in other words, the degree of the restriction of $\gamma_{\phi_t}$ to $\R \C_t$. Since $\chi(\tilde{V})=\chi(\C_t)$,  it follows from  Proposition \ref{degloggauss} that $\gamma_{\phi_t}$ is totally real. According to Proposition \ref{gaussmaptotreal}, the latter is equivalent to $\phi_t: \C_t \rightarrow \XD$ being a generalised  Harnack curve. The result follows.
\end{proof}

\section{Generalised Harnack curves with a single hyperbolic node}\label{sec:1hn}

\subsection{Tropical Harnack curves with a single hyperbolic node}

\begin{Proposition}\label{proptropharnsinglenode}
Let $h: C \rightarrow \mathbb{R}^2$ be a tropical Harnack curve with Newton polygon $\Delta$. Assume moreover that $h(C)$ has a single node $n$ that is hyperbolic, see Definition \ref{def:tropedgeleafnode}. Then, the parallelogram $n^\vee$ dual to $n$ in $\Sub_{h(C)}$ has exactly three of its vertices on the boundary of $\Delta$; these three vertices are distributed on two sides of $\Delta$ that intersect at smooth vertex $\vv$ of $\Delta$.
\end{Proposition}

\begin{proof}
Suppose to the contradiction that $n^\vee$ has at least two vertices $v_1$ and $v_2$ in the interior of $\Delta$. Consider the polygonal domain 
$P\subset\Delta$ obtained by taking the union of $n^\vee$ with all $2$-cells of $\Sub_{h(C)}$ 
having either $v_1$ or $v_2$ as a vertex. The only lattice points in the interior of $P$ are $v_1$ and $v_2$. It follows that the subset in $h(C)$ dual to $P$ has exactly two cycles. Consequently, the preimage of the latter subset in $C$ has only one cycle $\tilde{\lambda}$. There are two cases : either $v_1$ and $v_2$ are consecutive or opposite in $n^\vee$. In the first 
case, the set $\Gamma_\lambda$ (see Definition \ref{deftropsimpleHarnack}) consists of a single edge, namely the edge dual to the edge of $n^\vee$ joining $v_1$ to $v_2$. In the second case, the set $\Gamma_\lambda$ consists of the two edges 
crossing at the node $n$. In both cases, the condition of Definition \ref{deftropsimpleHarnack} is not fulfilled. This leads to a 
contradiction. It follows that $n^\vee$ has at most one vertex outside $\partial \Delta$.

The parallelogram $n^\vee$ cannot have four vertices on $\partial \Delta$, otherwise $C$ would be reducible and tropical Harnack curves are irreducible by definition. We deduce that $n^\vee$ has exactly three  vertices on $\partial \Delta$. Now, it follows from elementary observations that these three vertices are distributed on exactly two sides of $\Delta$ and that these two sides of $\Delta$ intersect at a smooth vertex of $\Delta$.
\end{proof}

\begin{Definition}\label{deftropnodenext}
Let $h:C \rightarrow\mathbb{R}^2$ be a tropical  Harnack curve with a single hyperbolic node $n$. We say that the node of $C$ is next to $\vv$ if  $\vv$ is the smooth vertex of the Newton polygon of $C$ given in Proposition \ref{proptropharnsinglenode}.
\end{Definition}

\begin{Proposition}\label{proptoptrop}
For any smooth vertex $\vv\in \Delta$ and any pair $h_1$, $h_2$ of tropical Harnack curves with Newton polygon $\Delta$ and with a single hyperbolic node next to $\vv$, we have 
\[\Top(h_1)=\Top(h_2).\]
We denote by $\Top(\Delta,\vv)$ the topological type $\Top(h)$ for any tropical Harnack curve $h$ with Newton polygon $\Delta$ and with a single hyperbolic node next to $\vv$.
\end{Proposition}

\begin{proof}
Let $h:C\rightarrow\R^2$ be a tropical Harnack curve with Newton polygon $\Delta$ and with a single hyperbolic node next to $\vv$. There is no loss of generality in assuming that $\vv=(0,0)$ and that the two sides of $\Delta$ adjacent to $\vv$ are supported on the positive coordinate axes. From Proposition \ref{proptropharnsinglenode}, we know that $n^\vee=\conv\big((0,1),(1,1),(k,0),(k+1,0)\big)$, up to permutation of the coordinates. In particular, the preimage $h^{-1}(n)$ cuts $C$ into a graph of genus $g_\Delta-1$ and two trees: one of these trees is the portion of an edge mapped vertically to $\R^2$ and the other one is a tree whose image in $\R^2$ is dual the the subdivided triangle $\conv\big((0,1),(0,0),(k,0)\big)$ in $\Sub_{h(C)}$. By making the bounded edges of the latter tree as small as desired, one can continuously deform the map $h$ so that the two trees of $C\setminus h^{-1}(n)$ have no vertices, that is $n^\vee=\conv\big((0,1),(1,1),(0,0),(0,0)\big)$. Along the deformation of $h$, the topological type $\Top(h)$ remains unchanged.
For any choice of sign, we can now reconstruct $\R V$ from $h$ as explained in Section \ref{subsecthc}. It is now clear that the topological type of $\R V$ does not depend on $h$, provided that $n^\vee=\conv\big((0,1),(1,1),(0,0),(0,0)\big)$.
\end{proof}

\subsection{Main statements}

In this section we undertake the classification of the topological pairs 
\[ \Top(\phi):= \Big( \R \TD, \; \R \C  \bigcup_{j\in \Z/n\Z} \R \mathcal{D}_j \Big)  \]
for generalised  Harnack curves $\phi:\C\rightarrow\XD$ with a single hyperbolic node. The main result is the following.

\begin{Theorem}\label{thmtopclass}
Let $\phi:\C\rightarrow\XD$ be a generalised Harnack curve with a single hyperbolic node. Assume moreover that $\C$ intersects transversally every toric divisor at infinity. Then, 
there is a unique smooth vertex $\vv$ of $\Delta$ such that 
\[\Top(\phi): = \Top(\Delta, \vv),\]
see Proposition \ref{proptropharnsinglenode} for the definition of \,$\Top(\Delta, \vv)$.
\end{Theorem}

\begin{Remark}
For any affine linear transformation $A$ sending $\Delta$ to itself, we have $\Top(\Delta, \vv)=\Top\big(\Delta, A(\vv)\big)$.
\end{Remark}

In the latter theorem, one had to specify the profile of intersection of $\phi(\C)$ with the toric divisors since the topological type $\Top(\phi)$ depends on it. Transversality is a generic condition and the general case can be deduced easily from the generic one.

\begin{figure}[h]
\centering
\input{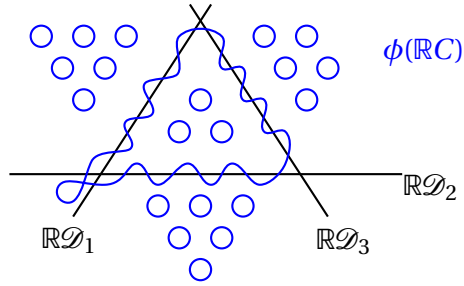}
\caption{A generalised Harnack curve of degree $8$ in $\cp{2}$ with a single hyperbolic node.}
\end{figure}

\begin{Definition}\label{defnodenext}
Let $\C \subset \TD$ be a generalised  Harnack curve with a single hyperbolic node. The node of $\C$ is said to be next to $\vv$ if  $\vv$ is the smooth vertex of $\Delta$ such that $\Top(\phi): = \Top(\Delta, \vv)$. 
\end{Definition}

\begin{Theorem}\label{thmgoodspine}
Let $\phi : \oC \rightarrow \XD$ be a generalised Harnack curve with a single hyperbolic node next to $\vv$. Then, there exist two abstract tropical curves $C_1, C_2 \subset \C$ exchanged by complex conjugation such that

-- $\Ap^{-1}\big(\mathscr{S}_{\phi(\C)}\big)=C_1\cup C_2$ where $\mathscr{S}_{\phi(\C)}$ is the spine of $\phi(\C)$,

-- both $\Ap: C_1 \rightarrow \R^2$ and $\Ap: C_2 \rightarrow \R^2$ are tropical Harnack curves with a single hyperbolic node next to $\vv$.
\end{Theorem}

\begin{Remark}

--It is not known in general if the spine of an algebraic curve of geometric genus $g$ in $\ttor$ can always be parametrised by an abstract tropical curve of genus $g$. Theorem \ref{thmgoodspine} implies that this is the case for generalised Harnack curves with a single hyperbolic node.

-- In \cite{KO}, it is shown that the map that associates the spine to a simple Harnack curve is a local diffeomorphism. If the same is true for generalised Harnack curves with a hyperbolic node, it would allow to study the deformations of such curves and answer the following questions. Can we contract the image of any compact connected component of $\R \C$ to a real isolated double point among the space of generalised Harnack curves with a hyperbolic node? Can we contract the self-intersecting loop of $\phi(\R \C)$ to a cusp among the same space of curves?
\end{Remark}

\subsection{Proofs}

In this section, the curve $\phi: \oC \rightarrow  \TD $ is a generalised Harnack curve with Newton polygon $\Delta$ and with a single hyperbolic node. We assume moreover that $\phi(\oC)$ intersects every toric divisor transversally, that is $\#\C_\infty=b_\Delta$. All the connected components $c_0,\dots,c_k$ of $\R\oC$ are oriented so that $\Ap(c_j)$ is locally convex.

\tocless\subsubsection{Topological maximality} 

\begin{Proposition}\label{proptopmax}
The set of points $\C_\infty$ is contained in a single connected component $c_0\subset\R\oC$. Moreover, the component $c_0$ contains the two branches of the hyperbolic node.
\end{Proposition}

\begin{figure}[h]
\centering
\input{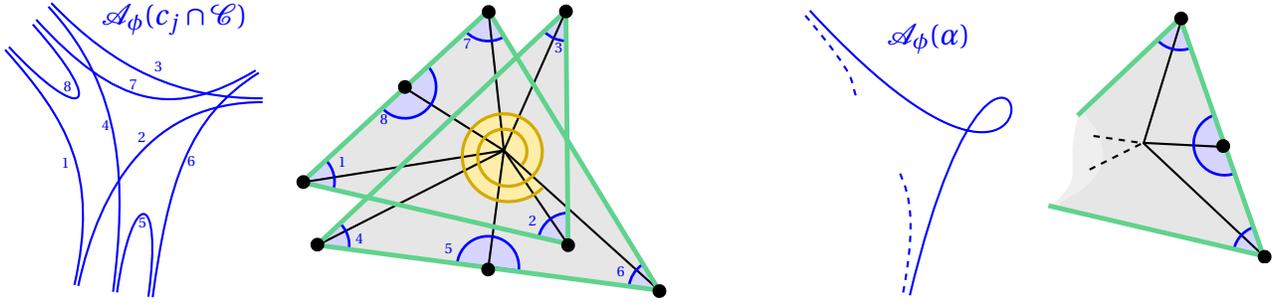}
\caption{Equation \eqref{eq:degcj} on the left and \eqref{eq:degalpha} on the right. The degree of $(\gp)_{\vert c_j}$ is the sum of the blue angles over $\pi$. The sum of the yellow angles is $2\pi r$ for a strictly positive integer $r$. The blue and yellow angles add up to $\pi$ times the number of green segments, that is $\pi\cdot\#(\C_\infty \cap c_j)$. In the presence of a component $\alpha$ as pictured on the left, we need to add $2$ segments.}
\label{fig:sumangle}
\end{figure}

\begin{proof}
In the proof of Theorem \ref{thm:max}, we constructed a deformation $\ell$ of $\R\oC$ within $\C_{\geq 0}$ whose connected components $\ell_0, \dots, \ell_k$ are in correspondence with the connected components $c_0,\dots,c_k$ of $\R\oC$ and such that the restriction of $\Ap$ to any $\ell_j$ is an immersion. Observe that if $c_j \cap \C_\infty=\emptyset$, then $c_j=\ell_j$ by construction. We claim the following:

-- if  $c_j \cap \C_\infty=\emptyset$, then $\ind\big(\Ap(c_j)\big)=1$ if $\phi(c_j)$ is smooth and $\ind\big(\Ap(c_j)\big)=2$ if $\phi(c_j)$ contains the two branches of the hyperbolic node,

-- if  $c_j \cap \C_\infty\neq\emptyset$, then $\ind\big(\Ap(\ell_j)\big)\leq 0$ with equality only if $c_j$ contains the two branches of the hyperbolic node. If moreover $\ind\big(\Ap(c_j)\big)<-1$, then $\phi(c_j)$ necessarily self-intersects.

Suppose to the contradiction that $c_0$ and $c_1$ have non-empty intersection with $\C_\infty$. Then, either the two branches of the node lift to a component $c_j$ disjoint from $\C_\infty$  and we deduce from the above claims that
\[ \ind\big(\Ap(\ell_0)\big)+\dots+\ind\big(\Ap(\ell_k)\big)\leq 2\cdot(-1)+(k-2)\cdot 1+1\cdot 2=k-2; \]
or the branches do not lift to a single such component $c_j$ in which case
\[ \ind\big(\Ap(\ell_0)\big)+\dots+\ind\big(\Ap(\ell_k)\big)\leq 1\cdot(-1)+1\cdot 0+(k-1)\cdot 1=k-2. \]
In both case, we obtain a contradiction with equation \eqref{eq:kauf}. The first part of the statement follows.

Let $c_0$ be the component containing $\C_\infty$. In order to satisfy equation \eqref{eq:kauf}, there are two possibilities left: either the two branches of the hyperbolic node lift to a single component $c_j$, $j\geq1$, or the branches lift to $c_0$. Indeed, the number of intersection point of any component $\phi(c_j)$, $j\geq1$,  with any other component $\phi(c_l)$ is necessarily even and $\phi(\oC)$ has only one node by assumption. In the first case, it follows from the first claim and equation \eqref{eq:kauf} that $\ind\big(\Ap(\ell_0)\big)=-2$. We deduce from the second claim that $\phi(c_0)$ has to self-intersect and consequently that the image of $\phi$ has at least two nodes. This is a contradiction. It follows that the two branches of the hyperbolic node lift to $c_0$. In particular, we have that $\ind\big(\Ap(\ell_j)\big)=1$ for $j\geq 1$ and then that $\ind\big(\Ap(\ell_0)\big)=-1$ in order to satisfy \eqref{eq:kauf}.

It remains to prove the two claims. For the first one, the parametrised curve $\Ap(c_j)$ has strictly positive curvature. In particular, the number of self-intersection of the curve is given by $\ind\big(\Ap(c_j)\big)-1$. Let us now prove the second claim. Assume that  $c_j$ is such that $\Ap$ is an embedding on each connected component of $c_j\setminus \C_\infty$. By summing the contribution to $\gp$ of each connected component of $c_j\setminus \C_\infty$, we obtain that 
\begin{equation}\label{eq:degcj}
\deg(\gp)_{\vert c_j} + 2\cdot r= \#(\C_{\infty}\cap c_j)
\end{equation}
where $r$ is a strictly positive integer, see Figure \ref{fig:sumangle}. According to equation \eqref{eq:anycc}, the index $\ind\big(\Ap(\ell_j)\big)$ is equal to $-r\leq -1$. Assume now that there is a connected component $\alpha \subset c_j\setminus \C_\infty$ such that $\Ap(\alpha)$ self-intersects once. We see, using similar arguments, that
\begin{equation}\label{eq:degalpha}
\deg(\gp)_{\vert c_j} + 2\cdot r = \#(\C_{\infty}\cap c_j)+2
\end{equation}
and deduce that $\ind\big(\Ap(\ell_j)\big)=1-r$, see again Figure \ref{fig:sumangle}.
Finally, if $\ind\big(\Ap(\ell_j)\big)\leq -2$, then there has to exist a connected components $\alpha$ of $c_j\setminus \C_\infty$ such that $\Ap(\alpha)$ intersects a tentacle of $\Ap(c_j\setminus \C_\infty)$ in two points. One of the two intersection points comes necessarily from a self-intersection of $\phi(c_j)$. This concludes the proof of the second claim as well as the proof of the statement.
\end{proof}

\begin{Corollary}\label{cor:ind}
Under the assumption of Proposition \ref{proptopmax}, we have that $\ind\big(\Ap(\ell_0)\big)=-1$. Moreover, the two branches of $\R\oC$ mapping to the hyperbolic node of $\phi(\oC)$ belong to a single connected component $\alpha \subset c_0\cap \C$.
\end{Corollary}

\begin{proof}
The equality $\ind\big(\Ap(\ell_0)\big)=-1$ is given in the core of the above proof. For the second part of the statement, assume to the contradiction that $\Ap$ is an embedding on each connected component of $c_0\cap \C$. Then, for two distinct components $\alpha, \beta \subset c_0\cap \C$, the convex arcs $\Ap(\alpha)$ and $\Ap(\beta)$ intersect in either $0$ or $2$ points. In the latter case, either none of the points corresponds to a node of $\phi(\R \C)$ or both do. In particular, the number of nodes of $\phi(\C)$ is even. This is a contradiction.
\end{proof}

\tocless\subsubsection{Lifted coamoebas}

Let $\vv$ be a smooth vertex of $\Delta$. Up to an affine change of coordinates, we can assume that $\vv=(0,0)$ and that the two edges of $\Delta$ adjacent to $\vv$ are supported on the positive horizontal and vertical axes. Denote by $\Delta_\vv \subset \Delta$ the polygonal domain defined by $$\Delta_\vv= \overline{\Delta \setminus \conv\big((0,0), (0,1),(1,0),(1,1)\big)}.$$ In the statement below, we denote by $\tau: \R^2\rightarrow\R^2$ the composition of the rotation by $-\pi/2$ and the homothety by $\pi$, both centred in the origin.

\begin{Proposition}\label{propliftcoa}
The restriction of the to argument map $\Arg$ to $\C_{>0}$ lifts to the universal covering $\R^2$ of $\sone$. Moreover, the closure of the lift of $\C_{>0}$ is a translation of $\tau \big( \Delta_\vv \big)$ for a unique smooth vertex $\vv \in \Delta$ .
\end{Proposition}

\begin{proof}
Let us first show that the restriction of $\Arg$ to $\C_{>0}$ lifts to the universal covering $\R^2$. To see this, it suffices to show that the induced map $\Arg_\ast: \pi_1(\C_{>0}) \rightarrow \pi_1(\sone)$ is trivial. According to Theorem \ref{thm:max} and Proposition \ref{proptopmax}, the set $\C_{>0}$ is an open disc with holes, each hole being bounded by a component of $\R\oC$ disjoint from $\C_\infty$. Since the map $\Arg$ is constant on each such component of  $\R\oC$, it follows that the map $\Arg_\ast$ is identically zero.

Let us now  determine the lift of $\C_{>0}$ to $\R^2$. Consider the deformation $\ell_0$ of the  component $c_0 \subset \R\oC$ containing $\C_\infty$ constructed in the proof of Theorem \ref{thm:max}. Recall that $\ell_0$ coincide with $c_0$ outside of an $\varepsilon$-neighbourhood of $\C_\infty$. When $\varepsilon\rightarrow 0$, the limit of the lift of $\ell_0$ to $\R^2$ is the image by $\tau$ of a piecewise linear curve $\Gamma\subset \R^2$ obtained by concatenation of all the primitive integer segments on $\partial \Delta$, according to Lemma \ref{extarg}. If follows also from the the latter lemma that the curve $\Gamma$  is the green curve of Figure \ref{fig:sumangle}, up to rescaling of each edge. In particular, the rotational index of $\Gamma$ is the same as the one of $\ell_0$ in absolute value, that is $1$ according to Corollary \ref{cor:ind}. Now, Corollary \ref{coralga} implies that $\Gamma$ has to enclose a domain of $\R^2$ which area, counted with multiplicities, is $\pi^2(\Area(\Delta)-1)$ (we applied a factor $1/2$ since $\C_{> 0}$ is only one half of $\C$). 

We claim that the latter implies that $\Gamma$ is the boundary of $\tau(\Delta_\vv)$ for a smooth vertex $\vv$ of $\Delta$. Indeed, for any pair of consecutive segments in $\Gamma$ that are not in convex position, the area of the polygonal domain enclosed by the curve obtained from $\Gamma$ be permuting these two segments increases at least by $\pi^2$ (if the segments happen to have opposite direction, we do two permutations at once). After finitely many such permutations, the piecewise linear curve we recover is $\tau(\partial \Delta)$. In particular, the curve $\tau(\partial \Delta)$ encloses a domain of area $\pi^2(\Area(\Delta)$. It follows that $\Gamma$ has to be obtained from $\tau(\partial \Delta)$ by a single permutation at a smooth vertex.
\end{proof}

\begin{Remark} The idea of lifting coamoebas to the universal cover of $\sone$ already appeared in \cite{PasNil}.
\end{Remark}

\begin{Corollary}\label{cor:loop}
Let $\vv \in \Delta$ be the smooth vertex given by Proposition \ref{propliftcoa}. Let $j \in \Z/n\Z$ such that $\vv=\Delta_j \cap \Delta_{j+1}$. Then, the (oriented) connected component $\alpha\subset c_0\cap \C$ of Corollary \ref{cor:ind} joins $\R \D_{j+1}$ to  $\R \D_{j}$.
\end{Corollary}

\begin{proof}
The connected component $\beta\subset c_0\cap C$ that maps to the non-convex lattice point of $\tau(\partial \Delta_\vv)$ joins $\R \D_{j+1}$ to  $\R \D_{j}$ when oriented so that $\Ap(\alpha)$ is convex. Hence, it has to self-intersect and $\alpha=\beta$.
\end{proof}

\begin{Corollary}\label{cor:order}
As above, let $\vv=\Delta_j \cap \Delta_{j+1}$ be the smooth vertex given by Proposition \ref{propliftcoa}. Then, the cyclic ordering on the set $\C_\infty$ induced by $c_0$ is as follows:  the oriented oval $c_0$ meets all the points of $\C_1$, then all the points of $\C_2$ and so on until $\C_n$, the only exception being that the last point of $\C_j$ encountered on $c_0$ lies between the two first points of $\C_{j+1}$.  
\end{Corollary}

\begin{proof}
All the components of $c_0\cap \C$ distinct from $\alpha$ embed in $\R^2$ as convex arcs. It follows that each such component has to join a point of $\C_l$ to a point of $\C_m$ such that the angle from the side $\Delta_l$ of $\Delta$ to the side $\Delta_m$ is strictly less than $\pi$. From the previous corollary, we know that $\alpha$ joins $\C_{j+1}$ to $\C_j$. We know from Corollary \ref{cor:ind} that $\ind\big(\A(\ell_0)\big)=-1$. The cyclic ordering of the statement is the only one to satisfy the latter conditions.
\end{proof}

\tocless\subsubsection{Spines}

Up to a toric transformation, one can assume that $\partial \Delta$ contains $(0,0)$, $(0,1)$ and $(1,0)$ and that the vertex $\vv$ given by Proposition \ref{propliftcoa} is $(0,0)$. Consider the affine chart $\mathbb{C}^2$ of $\TD$ whose coordinate axes are the two toric divisors adjacent to $\vv$. By Corollary \ref{cor:loop}, the arc $\A(\alpha)$ is a convex self-intersecting arc with a vertical asymptote at $x=-\infty$ and a horizontal asymptote at $y=-\infty$. Let $p=(p_1,p_2) \in \R^2$ be the point where $\A(\alpha)$ self-intersects and choose $\varepsilon_1, \, \varepsilon_2 >0$ such that the point $( p_1+ \varepsilon_1, p_2+ \varepsilon_2) $ belongs to the compact connected component of $\R^2 \setminus \mathcal{A}(\alpha)$. Define the following sets
\[
\begin{array}{c}
R := \left\lbrace (x,y) \in \mathbb{R}^2 \big| \: x \leq p_1 + \varepsilon_1, \: y \leq p_2 + \varepsilon_2 \right\rbrace, \\
H := \left\lbrace (x,y) \in \mathbb{R}^2 \big| \: x \leq  p_1   + \varepsilon_1, \: y =  p_2 + \varepsilon_2 \right\rbrace, \\
V := \left\lbrace (x,y) \in \mathbb{R}^2 \big| \: x =   p_1  + \varepsilon_1, \: y \leq  p_2   + \varepsilon_2 \right\rbrace.
\end{array}
\]
We choose $(\varepsilon_1,\varepsilon_2)$ generic such that both $H\cap \Ap(\C)$ and $V\cap \Ap(\C)$ are compact. 

\begin{Lemma}\label{lemma2to1}
The restriction of $\Ap$ to $\C \setminus \Ap^{-1}( R)$ is at most $2$-to-$1$. In particular, the restriction of $\Ap$ to $\R \C \setminus \Ap^{-1}(R)$ is an embedding, $\partial \big(\Ap(\C)\setminus R\big)= \Ap\big(\R \C \setminus \Ap^{-1}(R)\big)$ and $p$ belongs to $\R^2 \setminus \Ap(\C)$.
\end{Lemma}

\begin{proof}
Here, we will mainly follow the argumentation of \cite[Lemma 8]{Mikh}. The restriction of $\Ap$ to any connected component of $\R \C \setminus \alpha$ is an embedding since the restriction of $\phi$ is. The image of each such component by $\Ap$ cuts the plane into a convex and a non-convex half. The component $\Ap(\alpha)$ cuts the plane into three components. Only one of them is convex. It is moreover compact and contains $p$ in its interior. Among the two remaining non-convex component, one is contained in $R$, according to Corollary \ref{cor:loop}. It follows that the image by $\Ap$ of each component of $\R \C \cap \Ap^{-1}(\R^2 \setminus R)$ cuts $\R^2 \setminus R$ into a convex and a non-convex half.  Now, for any point $q\in \Ap(\C) \setminus R$ not in $\Ap(\R \C)$, draw a line segment $L\subset (\R^2 \setminus R)$ from $q$ to a point $r$ in a non-compact component of $\R^2 \setminus (\Ap(\C)\cup R)$. If $q$ belong to $a$ convex halves, then $\#\Ap^{-1}(q)=2(b-a)$ where $b$ is the number of convex halves containing $r$. From Corollary \ref{cor:order} and the construction of $R$, we deduce that $b=1$. We deduce that $\Ap$ is at most $2$-to-$1$ on $\C\setminus \Ap^{-1}(R)$.

For the second part of the statement, observe that $\Ap(\R \C)$ is the critical locus of $\Ap$, see Lemma \ref{lem:critloc}. Hence, the restriction of $\Ap$ to $\R \C\setminus \Ap^{-1}( R)$ is an embedding. The inclusion $\partial \big(\Ap(\C)\setminus R\big)\subset \Ap\big(\R \C \setminus \Ap^{-1}(R)\big)$ cannot be strict, otherwise there would be points in $\R^2 \setminus R$ with more than $2$ preimages by $\Ap$.  For $q=p$, we have that $a\geq1$ and $\#\Ap^{-1}(q)=2(1-a)$. It follows that $a=1$ and $\#\Ap^{-1}(q)=0$.
\end{proof}

\begin{Lemma}\label{lem:2to1}
One has that $\oC \cap \Ap^{-1}(R) = \C_H \cup \C_V$ where $ \C_H $ and $ \C_V$ are disjoint Riemann surfaces with boundary such that $ \Ap ( \partial \C_H ) \subset H $ and $ \Ap ( \partial \C_V ) \subset V $. Moreover, the restriction of $\Ap$ to either $ \C_H $ or $ \C_V$  is at most $2$-to-$1$.
\end{Lemma}

\begin{proof}
Recall that $V\cap\Ap(\C)$ is compact and $p$ is not in $\Ap(\C)$. Let $v_1$ (respectively $v_2$) be the uppermost (respectively lowermost) point on $V\cap \Ap(\C)$. By the previous lemma, both $v_1$ and $v_2$ sits on $\Ap(c_0)$. Consider any path $\rho_V$ joining $v_1$ to $v_2$ inside $\itr\big(\Ap(\C)\big) \setminus R$. Then, the lift $\Ap^{-1}(\rho_V)$ is a loop globally invariant by complex conjugation and cutting $c_0$ in exactly $2$ points, according to the previous lemma. By a symmetric construction, we obtain a path $\rho_H$ whose lift $\Ap^{-1}(\rho_H)$ is a loop globally invariant by complex conjugation and cutting $c_0$ in exactly $2$ points. The paths $\rho_V$ and $\rho_H$ can be chosen disjoint. In particular, the loops $\Ap^{-1}(\rho_H)$ and $\Ap^{-1}(\rho_V)$ cut $\oC$ in $3$ connected components: one contains $\Ap^{-1}(H)$, one contains $\Ap^{-1}(V)$ and the remaining component is disjoint from $\Ap^{-1}(R)$. If follows that we can decompose $\oC \cap \Ap^{-1}(R) = \C_H \cup \C_V$ as claimed in the statement.

For the second part of the statement, observe that there is a unique arc in both $\R \C_H$ and $\R \C_V$ that joins the two coordinate axes of the affine chart $\mathbb{C}^2$ (this is an easy consequence of Corollary \ref{cor:order}). We conclude by repeating the same arguments as in the proof of the previous lemma to the restriction of $\Ap$ to $\C_H$ and $\C_V$ respectively: for any point $r=(-R,-R)\in\R^2$ with sufficiently big $R>0$, we have that $r$ belongs to the complement of $\Ap(\C)$ and that $b=1$.
\end{proof}

\begin{Lemma}
There exist two holomorphic functions $g$ and  $h$  on $\A^{-1}(R)$ such that $\phi(\C_V)$ (respectively $\phi(\C_H)$) is the vanishing locus of $g$ (respectively of $h$) and $g \cdot h = f $ where $f$ is a defining polynomial for $\phi(\C)$.
\end{Lemma}

\begin{proof}
The closure of $\mathcal{A}^{-1}(R)$ in $\TD$ is a polydisc $D$ centred at the origin of the affine chart $\mathbb{C}^2$. In particular, the polydisc $D$ is a Stein manifold with $H^2(D,\mathbb{Z})=0$. According to \cite[Lemma 12, \S VIII.B]{GR}, there exist holomorphic functions $g$ and $h$ defined on $D$ whose vanishing locus are $\phi(\C_V)$ and $\phi(\C_H)$ respectively. Obviously, the functions $g$ and $h$ can be taken so that $f_{\vert D}=g\cdot h$.
\end{proof}

The two functions $g$ and $h$ admit respective Ronkin functions $N_g$ and $N_h$. Repeating the construction of Section \ref{sec:loggeo}, we obtain the spines $\mathscr{S}_g$ and $\mathscr{S}_h$ in $R$ from  piecewise linear function $S_g$ and $S_h$. By convexity of the Ronkin function, we have the inclusions $\mathscr{S}_g\subset \Ap(\C_V)$ and $\mathscr{S}_h\subset \Ap(\C_H)$ and the spines are deformation retract of the respective amoebas. The order maps of $g$ and $h$ of \cite{FPT} make sense and are still injective. The Legendre transform induces a duality between the subdivision of $R$ induced by $S_g$ (respectively $S_h$) and a subdivided lattice polygons $\Delta_g$ (respectively $\Delta_h$).

\begin{Lemma}\label{lemstabint}
The stable intersection of $\mathscr{S}_g$ and $\mathscr{S}_h$ in $R$ consists of a single point of multiplicity $1$.
\end{Lemma}

\begin{proof}
We can compose $\phi$ with an appropriate rescaling of the map $H_t$ of Section \ref{sec:ptcurve}, that we denote all together by $\phi_t$, such that $\phi_t(\C_V)$ and $\phi_t(\C_H)$ converge in Hausdorff distance to (unparametrised) phase-tropical curves in $\A^{-1}(R)$, see \cite[Proposition 6.1(3)]{Mikh05}. In the latter construction, we can require that the Hausdorff limit of $\A\big(\phi_t(\C_V)\big)$ and $\A\big(\phi_t(\C_H)\big)$ are $\mathscr{S}_g$ and $\mathscr{S}_h$ respectively. Along the way, we can guarantee that $ \phi_t(\C_V) \cap \A^{-1}(H) = \phi_t(\C_H) \cap \A^{-1}(V)=\emptyset$ so that the intersection multiplicity of $\phi_t(\C_V)$ with $\phi_t(\C_h)$ is $1$ for all $t$. In order to conclude, it suffices to show that the sum of the multiplicities of the stable intersection of $\mathscr{S}_g$ and $\mathscr{S}_h$ is the intersection multiplicity of $\phi(\C_V)$ with $\phi(\C_H)$. Up to a small translation, we can assume that $\mathscr{S}_g$ and $\mathscr{S}_h$ intersect transversally. The result follows now from \cite[Proposition 2.13]{BIMS} applied at every intersection point of $\mathscr{S}_g$ with $\mathscr{S}_h$.
\end{proof}

\begin{Lemma}
There are three possibilities for $\Delta_h$ and $\Delta_g$:

-- $\Delta_g$ is a vertical segment and $\Delta_h$ is a triangle contained in an horizontal strip of height $1$ with an horizontal side at the bottom.

-- $\Delta_h$ is an horizontal segment and $\Delta_g$ is a triangle contained in a vertical strip of width $1$ with a vertical side to the left, 

-- $\Delta_h$ is an horizontal segment and $\Delta_g$ is a vertical segment.\\In particular, the spine $\mathscr{S}_h$ (respectively $\mathscr{S}_g$) is a tree intersecting $H$ (respectively $V$) only once.
\end{Lemma}

\begin{proof}
Since $\C_H$ contains an arc of $\alpha$ whose image by $\phi$ intersection the horizontal coordinate axis of the affine chart $\mathbb{C}^2$, the polygon $\Delta_h$ has an horizontal side at its bottom. Similarly, the polygon $\Delta_g$ has a vertical side to its left. According to the previous lemma, the spines $\mathscr{S}_h$ and $\mathscr{S}_h$ intersect with multiplicity $1$, that is the mixed volume of $\Delta_h$ with $\Delta_g$ is $1$, see for instance \cite[Exercise 2(4)]{BIMS}. Recall also that $\Ap(\C_H)$ does not intersect $V$ and so neither do $\mathscr{S}_h$. In particular, $\Delta_h$ has no vertical side to the right. Similarly, $\Delta_g$ has no horizontal side at the top. The only possibilities left for $\Delta_g$ and $\Delta_h$ are the ones described in the statement.
\end{proof}

\begin{proof}[Proof of Theorem \ref{thmgoodspine}] 
Since the spine $\mathscr{S}_h$ (respectively $\mathscr{S}_g$) is a tree intersecting $H$ (respectively $V$) only once, we deduce that $\R \C_H \subset c_0$ and $\R \C_V \subset c_0$. It implies that $\Ap$ maps the $g_\Delta-1$ compact components of $\R \C$ in the complement of $R$. Lemma \ref{lemma2to1} implies that the latter $g_\Delta-1$ ovals bound a compact component of $\R^2 \setminus \Ap(\C)$ and then, that the intersection of the spine of $\phi(\C)$ with the complement of $R$ has $g_\Delta-1$ cycles. We claim that the intersection of the spine of $\phi(\C)$ with $R$ is the union of $\mathscr{S}_h$ and $\mathscr{S}_g$. All together, it implies that $\Ap^{-1}\big(\mathscr{S}_{\phi(\C)}\big)$ is the disjoint union of $2$ tropical curves $C_1$ and $C_2$ exchange by complex conjugation, according to the $2$-to-$1$ property of Lemmas \ref{lemma2to1} and \ref{lem:2to1}. Since the preimage of the node $n$ of $\mathscr{S}_{\phi(\C)}$ is disjoint from the cycles of $C_1$ and $C_2$, both $\Ap: C_1 \rightarrow \R^2$ and $\Ap: C_2 \rightarrow \R^2$ are Harnack. The fact that $n$ is next to $\vv=(0,0)$ follows from Definition \ref{deftropnodenext} and the previous lemma. 

It remains to prove the claim. Since $\mathscr{S}_h$ and $\mathscr{S}_g$ are both trees, no component of the complement of $\Ap(\C_H)$ is hidden by $\Ap(\C_V)$ and vice versa. It follows that the piecewise linear function $S_f$ defining the spine $\mathscr{S}_{\phi(\C)}$ in $R$ is the sum of $S_g$ and $S_h$. It implies that $\mathscr{S}_{\phi(\C)}\cap R=\mathscr{S}_h\cup\mathscr{S}_g$.
\end{proof}

\begin{proof}[Proof of Theorem \ref{thmtopclass}]
Recall from Theorem \ref{thm:max} that the quotient of $\oC$ by its anti-holomorphic involution $\sigma$ is a disc with $g_\Delta-1$ holes, bounded by $\R \oC$. By construction,  the graphs $C_1$ and $C_2$ of Theorem \ref{thmgoodspine} are exchanged by $\sigma$ and project down to a skeleton of the quotient $\C / \sigma$. If we denote by $C$ the projection of $C_1$ to $\C / \sigma$, the inclusion $C \subset \C / \sigma$ is a ribbon graph such that the Harnack curve $\Ap: C \rightarrow \R^2$ extends into an immersion $\Ap: \C / \sigma \rightarrow \R^2$. The topological type of $\Top(\phi)$ can then be recovered from the procedure described in Section  \ref{subsecthc}. Since $\Ap: C \rightarrow \R^2$ is a tropical Harnack curve with a single hyperbolic node next to $\vv$ according to Theorem \ref{thmgoodspine}, it follows that 
\[\Top(\phi)= \Top\big(\Ap: C \rightarrow \R^2\big)=\Top(\Delta,\vv).\]
\end{proof}

\bibliographystyle{alpha}
\bibliography{Draft2}

\noindent Matematiska institutionen, Stockholms universitet, 106 91 Stockholm, Sweden\\
Email : lang@math.su.se

\end{document}